\documentclass[11pt]{article}
\usepackage{amsfonts,mathrsfs,amssymb,amsthm,mathptm}
\usepackage{amsmath,amscd}
\usepackage{color}
\usepackage[dvips]{graphicx}
\usepackage[all]{xy}
\usepackage{graphicx}
\usepackage{fancyhdr}
\usepackage[toc,page]{appendix}
\usepackage{MnSymbol}

\begin{document}
\newtheorem{Def}{Definition}[section]
\newtheorem{Bsp}[Def]{Example}
\newtheorem{Prop}[Def]{Proposition}
\newtheorem{Theo}[Def]{Theorem}
\newtheorem{Lem}[Def]{Lemma}
\newtheorem{Koro}[Def]{Corollary}
\theoremstyle{definition}
\newtheorem{Rem}[Def]{Remark}

\newcommand{\add}{{\rm add}}
\newcommand{\gd}{{\rm gl.dim}}
\newcommand{\dm}{{\rm dom.dim}}
\newcommand{\E}{{\rm E}}
\newcommand{\Mor}{{\rm Morph}}
\newcommand{\End}{{\rm End}}
\newcommand{\ind}{{\rm ind}}
\newcommand{\rsd}{{\rm resdim}}
\newcommand{\rd} {{\rm repdim}}
\newcommand{\ol}{\overline}
\newcommand{\overpr}{$\hfill\square$}
\newcommand{\rad}{{\rm rad}}
\newcommand{\soc}{{\rm soc}}
\renewcommand{\top}{{\rm top}}
\newcommand{\pd}{{\rm projdim}}
\newcommand{\id}{{\rm injdim}}
\newcommand{\fld}{{\rm flatdim}}
\newcommand{\Fac}{{\rm Fac}}
\newcommand{\Gen}{{\rm Gen}}
\newcommand{\Supp}{{\rm Supp}}
\newcommand{\Ass}{{\rm Ass}}
\newcommand{\Spec}{{\rm Spec}}

\newcommand{\fd} {{\rm fin.dim}}
\newcommand{\DTr}{{\rm DTr}}
\newcommand{\cpx}[1]{#1^{\bullet}}
\newcommand{\D}[1]{{\mathscr D}(#1)}
\newcommand{\Dz}[1]{{\mathscr D}^+(#1)}
\newcommand{\Df}[1]{{\mathscr D}^-(#1)}
\newcommand{\Db}[1]{{\mathscr D}^b(#1)}
\newcommand{\Ds}[1]{{\mathscr D}^{\ast}(#1)}
\newcommand{\C}[1]{{\mathscr C}(#1)}
\newcommand{\Cz}[1]{{\mathscr C}^+(#1)}
\newcommand{\Cf}[1]{{\mathscr C}^-(#1)}
\newcommand{\Cb}[1]{{\mathscr C}^b(#1)}
\newcommand{\K}[1]{{\mathscr K}(#1)}
\newcommand{\Kz}[1]{{\mathscr K}^+(#1)}
\newcommand{\Kf}[1]{{\mathscr  K}^-(#1)}
\newcommand{\Kb}[1]{{\mathscr K}^b(#1)}
\newcommand{\modcat}{\ensuremath{\mbox{{\rm -mod}}}}
\newcommand{\Modcat}{\ensuremath{\mbox{{\rm -Mod}}}}

\newcommand{\stmodcat}[1]{#1\mbox{{\rm -{\underline{mod}}}}}
\newcommand{\pmodcat}[1]{#1\mbox{{\rm -proj}}}
\newcommand{\imodcat}[1]{#1\mbox{{\rm -inj}}}
\newcommand{\Pmodcat}[1]{#1\mbox{{\rm -Proj}}}
\newcommand{\Imodcat}[1]{#1\mbox{{\rm -Inj}}}
\newcommand{\opp}{^{\rm op}}
\newcommand{\otimesL}{\otimes^{\rm\mathbb L}}
\newcommand{\rHom}{{\rm\mathbb R}{\rm Hom}\,}
\newcommand{\projdim}{\pd}
\newcommand{\Hom}{{\rm Hom}}
\newcommand{\Coker}{{\rm Coker}}
\newcommand{ \Ker  }{{\rm Ker}}
\newcommand{ \Img  }{{\rm Im}}
\newcommand{\Ext}{{\rm Ext}}
\newcommand{\StHom}{{\rm \underline{Hom}}}

\newcommand{\gm}{{\rm _{\Gamma_M}}}
\newcommand{\gmr}{{\rm _{\Gamma_M^R}}}

\def\vez{\varepsilon}\def\bz{\bigoplus}  \def\sz {\oplus}
\def\epa{\xrightarrow} \def\inja{\hookrightarrow}

\newcommand{\ra}{\rightarrow}
\newcommand{\lra}{\longrightarrow}
\newcommand{\lraf}[1]{\stackrel{#1}{\lra}}
\newcommand{\lla}{\longleftarrow}
\newcommand{\llaf}[1]{\stackrel{#1}{\lla}}
\newcommand{\dk}{{\rm dim_{_{k}}}}

\newcommand{\colim}{{\rm colim\, }}
\newcommand{\limt}{{\rm lim\, }}
\newcommand{\Add}{{\rm Add }}
\newcommand{\Prod}{{\rm Prod }}
\newcommand{\Tor}{{\rm Tor}}
\newcommand{\Cogen}{{\rm Cogen}}

{\Large \bf
\begin{center}
Derived decompositions of abelian categories, I.
\end{center}}
\medskip
\centerline{\textbf{Hong Xing Chen} and \textbf{Chang Chang Xi}$^*$}

\renewcommand{\thefootnote}{\alph{footnote}}
\setcounter{footnote}{-1} \footnote{ $^*$ Corresponding author.
Email: xicc@cnu.edu.cn; Fax: +86 10 68903637.}
\renewcommand{\thefootnote}{\alph{footnote}}
\setcounter{footnote}{-1} 
\renewcommand{\thefootnote}{\alph{footnote}}
\setcounter{footnote}{-1} 

\begin{abstract}
Derived  decompositions of abelian categories are introduced in internal terms of abelian subcategories to construct semi-orthogonal decompositions (or Bousfield localizations, or hereditary torsion pairs) in various derived categories of abelian categories. We give a sufficient condition for arbitrary abelian categories to have such derived decompositions and show that it is also necessary for abelian categories with enough projectives and injectives. For bounded derived categories, we describe which semi-orthogonal decompositions are determined by derived decompositions. The necessary and sufficient condition is then applied to the module categories of rings: localizing subcategories, homological ring epimorphisms, commutative noetherian rings and nonsingular rings. Moreover, for a commutative noetherian ring of Krull dimension at most $1$, a derived stratification of its module category is established.
\end{abstract}
\medskip
{\footnotesize\tableofcontents\label{contents}}

\section{Introduction}
Semi-orthogonal decompositions (or hereditary torsion pairs in the terminology of \cite{BI}) have been applied in a number of branches of mathematics. For example, in homotopy and triangulated categories, they were also named as Bousfield localizations (\cite[Section 9.1]{neemanbook}) and applied to get $t$-structures of triangulated categories (see \cite{BBD}), 
and in algebraic geometry they were used to study Fourier-Mukai transforms on derived categories of coherent sheaves of smooth projective varieties
(see \cite[Chapter 11]{Huybrechts}, \cite{Bondal}).
However, in the course of studying semi-orthogonal decompositions in triangulated categories, the following fundamental question seems to remain:

\smallskip
{\bf Question}. Given an abelian category $\mathcal{A}$, how can we construct semi-orthogonal decompositions (or hereditary torsion pairs ) of the $\ast$-bounded derived category $\Ds{\mathcal{A}}$ of $\mathcal{A}$ for $*\in \{b, +,-, \varnothing\}?$

\smallskip
By definition, semi-orthogonal decompositions of $\Ds{\mathcal{A}}$ are defined at the level of derived categories (see Definition \ref{SOD}), but we would like instead to have a characterization of such decompositions directly at the level of given abelian categories themselves. So we introduce naturally the notion of derived decompositions of abelian categories.

\begin{Def}\label{Intro-def}
Let $\mathcal{A}$ be an abelian category, and let $\mathcal{X}$ and $\mathcal{Y}$ be full subcategories of $\mathcal{A}$.  For $\ast\in \{b,+,-, \varnothing\}$, $\Ds{\mathcal{A}}$ denotes the $\ast$-bounded derived category of $\mathcal{A}$. The pair $(\mathcal{X}, \mathcal{Y})$ is called a \emph{$\mathscr{D}^*$-decomposition} of $\mathcal{A}$ if

$(D1)$ both $\mathcal{X}$ and $\mathcal{Y}$ are abelian subcategories of $\mathcal{A}$, and the inclusions $\mathcal{X}\subseteq\mathcal{A}$ and $\mathcal{Y}\subseteq\mathcal{A}$ induce fully faithful functors $\Ds{\mathcal{X}}\to \Ds{\mathcal{A}}$ and $\Ds{\mathcal{Y}}\to\Ds{\mathcal{A}}$, respectively.

$(D2)$ $\Hom_{\Ds{\mathcal{A}}}(X, Y[n])=0$ for all $X\in \mathcal{X}$, $Y\in\mathcal{Y}$ and $n\in\mathbb{Z}.$

$(D3)$ For any object $\cpx{M}\in\Ds{\mathcal{A}}$, there is a triangle in $\Ds{\mathcal{A}}$
$$
X_{\cpx{M}}\lra \cpx{M}\lra Y^{\cpx{M}}\lra X_{\cpx{M}}[1]
$$
such that $X_{\cpx{M}}\in\Ds{\mathcal{X}}$ and $Y^{\cpx{M}}\in\Ds{\mathcal{Y}}$.
\end{Def}

For convenience, a $\mathscr{D}^b$-decomposition of $\mathcal{A}$ is also termed \emph{derived decomposition} of $\mathcal{A}$ in the sequel.

In this paper we establish a characterization of $\mathscr{D}^*$-decompositions (thus also semi-orthogonal decompositions) in entirely internal terms of conditions on subcategories of given abelian categories, instead of the ones of derived categories (see Definition \ref{SOD}). The characterization is then applied explicitly to a wide variety of situations for module categories, including homological ring epimorphisms, localizing subcategories and commutative noetherian rings.

These applications motivate us to introduce derived stratifications of abelian categories (see Section \ref{sect2.3} for definition). We show that, among others, the module category of a commutative noetherian ring of Krull dimension at most $1$ has a derived stratification with abelian simple factors. But for an indecomposable commutative ring its bounded derived category does not have non-trivial stratification by bounded derived categories of rings (see \cite{AKLY}). Compared with this phenomenon, the notion of derived decompositions may be of interest for stratifying bounded derived categories of rings by bounded derived categories of abelian categories. This provides a way to approach the derived category of an abelian category by those of its smaller abelian subcategories.

Our main result reads as follows.

\begin{Theo}\label{main-result}
Let $\mathcal{A}$ be an abelian category, $\mathcal{X}$ and $\mathcal{Y}$ full subcategories of $\mathcal{A}$ and $\ast\in \{b,+,-, \varnothing\}$.

$(1)$ The pair $(\mathcal{X},\mathcal{Y})$ is a $\mathscr{D}^*$-decomposition of $\mathcal{A}$
if the following conditions hold:

\quad $(a)$ $\Ext_\mathcal{A}^n(X,Y)=0$ for any $n\geq 0$, $X\in\mathcal{X}$ and $Y\in\mathcal{Y}$.

\quad $(b)$ For each object $M\in\mathcal{A}$, there is a long exact sequence
$$0\lra Y_M\lra X_M\lra M\lra Y^M\lra X^M\lra 0$$
in $\mathcal{A}$ with $X_M, X^M\in\mathcal{X}$ and $Y_M, Y^M\in\mathcal{Y}.$

\quad $(c)$ For each object $M\in\mathcal{A}$, there is a monomorphism $M\to I$ in $\mathcal{A}$ such that $X^{I}=0$ in $(b).$

\quad $(d)$ For each object $M\in\mathcal{A}$, there is an epimorphism $P\to M$ in $\mathcal{A}$ such that $Y_{P}=0$ in $(b).$

\smallskip
$(2)$ Suppose that $\mathcal{A}$ has enough projectives and injectives. Then the pair $(\mathcal{X},\mathcal{Y})$ is a $\mathscr{D}^*$-decomposition of $\mathcal{A}$ if and only if the above $(a)$ and $(b)$ together with $(c')$ and $(d')$ hold, where

\quad $(c')$ if $M$ is injective, then $X^M=0$ in $(b).$

\quad $(d')$ If $M$ is projective, then $Y_M=0$ in $(b).$
\end{Theo}

Be aware that $(b)$ in Theorem \ref{main-result}(1) was introduced in \cite{Krause-S} to study the telescope conjecture for
hereditary rings. The two conditions (a) and (b) are equivalent to saying that $(\mathcal{X},\mathcal{Y})$ is a complete Ext-orthogonal pair in $\mathcal{A}$. Unfortunately, (b) has not been used elsewhere in the literature before our investigation in this paper.

For bounded derived categories, we characterize which semi-orthogonal decompositions are induced from derived decompositions (see Proposition \ref{Corre} for details).

Theorem \ref{main-result}(2) implies that if $\mathcal{A}$ has enough projectives and injectives,
then the existences of $\mathscr{D}^*$-decompositions of $\mathcal{A}$ for all $\ast\in \{b,+,-, \varnothing\}$ are equivalent.
This applies to the module categories of rings. In particular, we have the following
consequence of Theorem \ref{main-result} on homological ring epimorphisms.

\begin{Koro}
Let $\lambda: R\to S$ be a homological ring epimorphism. Define
$\mathcal{Y}:=\{Y\in R\Modcat\mid\Hom_R(S, Y)=0=\Ext_R^1(S, Y)\}$ and $\mathcal{Z}:=\{X\in R\Modcat\mid S\otimes_RX=0=\Tor_1^R(S, X)\}$. Then

$(1)$ $(S\Modcat, \mathcal{Y})$ is a derived decomposition of  $R\Modcat$ if and only if  $\pd(_RS)\leq 1$ and
$\Hom_R(\Coker(\lambda)$, $\Ker(\lambda))=0$.

$(2)$ $(\mathcal{Z}, S\Modcat)$ is a derived decomposition of $R\Modcat$ if and only if $\fld(S_R)\leq 1$ and $\Coker(\lambda)\otimes_RI=0$ for any injective $R$-module $I$.

$(3)$ If the conditions  $(1)$ and $(2)$ are satisfied, then, for any $\ast\in \{b,+,-, \varnothing\}$, $\Ds{\mathcal{Y}}\lraf{\simeq} \Ds{\mathcal{Z}}$ and there exists a recollement:
$$\xymatrix@C=1.3cm{\Ds{S}\ar[r]^-{\Ds{\lambda_*}}&\Ds{R}\ar[r]
\ar@/^1.2pc/[l]\ar@/_1.8pc/[l]
&\Ds{\mathcal{Y}}\ar@/^1.2pc/[l]\ar@/_1.8pc/[l]^-{}}$$
\end{Koro}

Applying Theorem \ref{main-result} to commutative rings, we have the following corollary. For notation and notions, we refer the reader to Section \ref{sect4.3}.

\begin{Koro}\label{cor1.4}
Let $R$ be a commutative noetherian ring.

$(1)$ Suppose that $\Phi$ is a specialization closed subset of $\Spec(R)$. If the Krull dimension of $R$ is at most $1$, then $\big(\Supp^{-1}(\Phi), \Supp^{-1}(\Phi^c)\big)$ is a derived decomposition of $R\Modcat$, where $\Phi^{c}:=\Spec(R)\setminus\Phi$.

$(2)$ Let $\Sigma$  be a multiplicative subset of $R$,  $\Sigma^-R$ the localization of $R$ at $\Sigma$ and $\Phi:=\{\mathfrak{p}\in\Spec(R)\mid\mathfrak{p}\cap\Sigma\neq \emptyset\}.$ Then
$\big(\Supp^{-1}(\Phi), (\Sigma^-R)\Modcat\big)$ is a derived decomposition of $R\Modcat$.
\end{Koro}

Further applications of Theorem \ref{main-result} to localizing subcategories and nonsingular rings, are given by Proposition \ref{Localizing} and Corollary \ref{DADS}, respectively. Note that, for a commutative noetherian ring of Krull dimension at most $1$, we show that a derived stratification by derived categories of abelian categories always exists (see Corollary \ref{decomposition}).

The article is outlined as follows: In Section \ref{sect2} we fix notation and recall definitions needed in proofs. In Section \ref{sect3} we prove Theorem \ref{main-result}. The proof is divided into two parts. The first one is for the proof of Theorem \ref{main-result}(1), while the second one is for that of Theorem \ref{main-result}(2). In Section \ref{sect4} we apply Theorem \ref{main-result} to construct derived decompositions of module categories from various aspects: ring epimorphisms, localizing subcategories, commutative noetherian rings and nonsingular rings. Particularly, we show the strong conclusion, Corollary \ref{decomposition}.

In the second paper we shall give a series of applications of derived decompositions. In particular, we construct complete cotorsion pairs from derived decompositions, with applications to infinitely generated tilting modules.

\medskip
{\bf Acknowledgement.} The both authors would like to thank Ch. Psaroudakis,  L. Positselski and M. Saorin
for helpful discussions, and NSFC for partial support.

\section{Notation and definitions\label{sect2}}
In this section we first fix some notation and recall definitions of semi-orthogonal decompositions (or hereditary torsion pairs),
cotorsion pairs and complete Ext-orthogonal pairs.

\subsection{Notation for derived categories}

Let $\mathcal A$ be an additive category.

A full subcategory $\mathcal B$ of
$\mathcal A$ is always assumed to be closed under isomorphisms.
For an object $X\in\mathcal{A}$,  $\add(X)$
(respectively, $\Add(X)$) denotes the full subcategory of
$\mathcal{A}$ consisting of all direct summands of finite
(respectively, arbitrary) coproducts of copies of $X$ (if arbitrary coproducts exist).

Let $F:\mathcal {A}\to \mathcal{A}'$ be an additive functor from $\mathcal{A}$ to another additive category $\mathcal{A}'$. The kernel and image of $F$ are defined as $\Ker(F):=\{X\in \mathcal{A}\mid FX\simeq 0\}$ and
$\Img(F):=\{Y\in \mathcal{A}'\mid \exists\, X\in \mathcal{A},\; FX\simeq Y\}$, respectively. Let $f:X\to Y$ be a morphism in $\mathcal{A}$. The kernel, image and cokernel of $f$, whenever they exist,
will be denoted by $\Ker(f)$, $\Img(f)$ and $\Coker(f)$, respectively.

By a complex $\cpx{X}$ over $\mathcal{A}$ we mean a sequence of
morphisms $d^i$ between objects $X^i$ in $\mathcal{A}:\; \cdots\to
X^i\epa{d^i } X^{i+1}\epa{d^{i+1}} X^{i+2}\to\cdots$, such that
$d^id^{i+1}=0$ for all $i\in\mathbb{Z}$. We write
$\cpx{X}=(X^i, d^i)_{i\in\mathbb{Z}}$ and call $d^i$ the $i$-th
differential of $\cpx{X}$. For a fixed $n\in \mathbb{Z}$, we denote by
$\cpx{X}[n]$ the complex obtained from $\cpx{X}$ by shifting $n$
degrees, that is, $(\cpx{X}[n])^i=X^{n+i}$ with the $i$-th differential $(-1)^nd^{n+i}$, and by $H^n(\cpx{X})$
the $n$-th cohomology of $\cpx{X}$.

Let $\C{\mathcal{A}}$ be the category of all complexes over
$\mathcal{A}$ with chain maps as morphisms, and $\K{\mathcal{A}}$ the homotopy
category of $\C{\mathcal{A}}$. We denote by $\Cb{\mathcal{A}}$ and
$\Kb{\mathcal{A}}$ the bounded complex and homotopy categoires of $\mathcal{A}$, respectively.

From now on, let $\mathcal{A}$ be an abelian category.

By $\D{\mathcal{A}}$ and $\Db{\mathcal{A}}$ we denote the \emph{unbounded} and \emph{bounded derived categories} of $\mathcal{A}$, respectively.
Throughout the paper, we always identify $\Db{\mathcal{A}}$ with the full subcategory of $\D{\mathcal{A}}$ consisting of all complexes with finitely many non-zero cohomologies because they are equivalent as triangulated categories. Further, by $\Dz{\mathcal{A}}$ and $\Df{\mathcal{A}}$ we denote the \emph{bounded-below} and \emph{bounded-above derived categories} of $\mathcal{A}$, respectively.

For any $X, Y\in\mathcal{A}$ and $i\in\mathbb{Z}$, we write $\Ext_\mathcal{A}^i(X,Y)$ for $\Hom_{\D{\mathcal{A}}}(X, Y[i])$. Note that $\Ext_\mathcal{A}^0(X,Y)=\Hom_\mathcal{A}(X,Y)$ and  $\Ext_\mathcal{A}^i(X,Y)=0$ whenever $i<0$. For each $i\geq 1$, $\Ext_\mathcal{A}^i(X,Y)$ can be identified with the set of equivalence classes of long exact sequences
$0\to Y\to E_i\to \cdots\to E_1\to X\to 0$ in $\mathcal{A}$ (see \cite[XI]{Iversen} for details).

The following facts are standard in homological algebra.

$(1)$ Suppose that $\mathcal{A}$ has enough projectives with $\mathscr{P}(\mathcal{A})$ the category of all projective objects of $\mathcal{A}$. Further, let $\mathscr{K}^{-, b}(\mathscr{P}(\mathcal{A}))$ be the full subcategory of $\K{\mathcal{A}}$ consisting of bounded-above complexes with all terms in $\mathscr{P}{(\mathcal{A})}$ and finitely many nonzero cohomologies. Then there is a triangle equivalence between
$\mathscr{K}^{-, b}(\mathscr{P}(\mathcal{A}))$ and $\Db{\mathcal{A}}$. In this case, $\Ext_\mathcal{A}^i(X,Y)$ is isomorphic to the usual $i$-th extension group of $X$ and $Y$, defined by projective resolutions of $X$.

$(2)$ Dually, suppose that $\mathcal{A}$ has enough injectives with $\mathscr{I}(\mathcal{A})$ the category of all injective objects of $\mathcal{A}$. Then there is a triangle equivalence between $\mathscr{K}^{+, b}(\mathscr{I}(\mathcal{A}))$ and $\Db{\mathcal{A}}$, where $\mathscr{K}^{+, b}(\mathscr{I}(\mathcal{A}))$ is defined similarly. In this situation, $\Ext_\mathcal{A}^i(X,Y)$ can be calculated by taking injective resolutions of $Y$.

\smallskip

A full subcategory $\mathcal{B}$ of $\mathcal{A}$ is called an \emph{abelian subcategory} of $\mathcal{A}$ if $\mathcal{B}$ is an abelian category and the inclusion $\mathcal{B}\to\mathcal{A}$ is an exact functor between abelian categories. This is equivalent to saying that $\mathcal{B}$ is closed under taking
kernels and cokernels in $\mathcal{A}$. The full subcategories $\{0\}$ and $\mathcal{A}$ are called the \emph{trivial} abelian subcategories of $\mathcal{A}$.

For $n\in\mathbb{N}$ and a full subcategory $\mathcal{B}$ of $\mathcal{A}$, we define the full subcategories of $\mathcal{A}$:

$$\begin{array}{rl} {^{\bot n}}\mathcal{B}:= & \{X\in \mathcal{A}\mid \Ext_{\mathcal{A}}^n(X,Y)=0, Y\in\mathcal{B}\}, \\
{^{\bot {>n}}}\mathcal{B} := & \{X\in \mathcal{A}\mid \Ext_{\mathcal{A}}^j(X,Y)=0, Y\in\mathcal{B}, j>n\},\\
{^\bot}\mathcal{B} := & \{X\in \mathcal{A}\mid \Ext_{\mathcal{A}}^j(X,Y)=0, Y\in\mathcal{B}, j\ge 0 \}.
\end{array}
$$
Similarly, $\mathcal{B}{^{\bot n}}$, $\mathcal{B}{^{\bot {>n}}}$ and  $\mathcal{B}^{\bot}$ are defined.
Recall that ${^\bot}\mathcal{B}$ is said to be \emph{left perpendicular} to $\mathcal{B}$ in $\mathcal{A}$, while
$\mathcal{B}{^\bot}$ is said to be \emph{right perpendicular} to $\mathcal{B}$ in $\mathcal{A}$ (see \cite{GL}).

Let $F:\mathcal{A}\to \mathcal{A}'$ be an exact functor of abelian categories. Then $F$ induces  derived functors $\Ds{F}:\Ds{\mathcal{A}}\to \Ds{\mathcal{A}'}$
for any $\ast\in \{b,+,-, \varnothing\}$, defined by $F(\cpx{X}):=(FX^i, Fd^i)_{i\in \mathbb{Z}}$ for $\cpx{X}\in\Ds{\mathcal{A}}$.

\begin{Lem}\label{subcat}
If $\mathcal{B}$ is a abelian subcategory of an abelian category $\mathcal A$ such that the inclusion $i: \mathcal{B}\ra \mathcal{A}$ induces a fully faithful functor $\Db{i}: \Db{\mathcal{B}}\ra \Db{\mathcal{A}}$, then $\Img(\Db{i})$ consists of all bounded complexes $\cpx{X}$ with $H^n(\cpx{X})\in\mathcal{B}$  for all $n\in\mathbb{Z}$.
\end{Lem}

{\it Proof.} Let $\mathscr{X}$ be the full subcategory of $\Db{\mathcal{A}}$ consisting of complexes with all cohomologies in $\mathcal{B}$. Note that $\Img(\Db{i})$ is a full triangulated subcategory of $\Db{\mathcal{A}}$ containing $\mathcal{B}$ and being closed under isomorphisms. Since $\mathcal{B}\subseteq \mathcal{A}$ is an ableian subcategory, $\Img(\Db{i})\subseteq \mathscr{X}$. Now, let $\cpx{X}\in \mathscr{X}$. We can show by induction on the number of nonzero terms of $\cpx{X}$ that $\cpx{X}$ belongs to the smallest full triangulated subcategory containing $H^n(\cpx{X})$ for all $n\in \mathbb{Z}$. Clearly, this subcategory is contained in $\Img(\Db{i})$. Thus $\cpx{X}\in \Img(\Db{i})$ and $\mathscr{X}\subseteq \Img(\Db{i})$. Hence $\Img(\Db{i})=\mathscr{X}.$ $\square$

\medskip
By a ring we mean an associative ring $R$ with identity. We denote by $R\Modcat$ the category of all
unitary left $R$-modules. For an $R$-module $M$, we denote by $\pd(_RM)$, $\id(_RM)$ and $\fld(_RM)$
the projective, injective and flat dimensions of $M$, respectively.
As usual, we simply write
$\C{R}$, $\K{R}$ and $\D{R}$ for the complex, homotopy and derived categories of $R\Modcat$, respectively.

Let $\lambda: R\ra S$ be a homomorphism of rings. We denote by $\lambda_*:S\Modcat\to R\Modcat$ the \emph{restriction functor} induced by $\lambda$, and by $\D{\lambda_*}:\D{S}\to\D{R}$ the derived functor of $\lambda_*$. If $\lambda_*$ is fully faithful, then $\lambda$ is called a \emph{ring epimorphism}.
If $\D{\lambda_*}$ is fully faithful, then $\lambda$ is called a \emph{homological} ring epimorphism.
Note that $\lambda$ is a homological ring epimorphism if and only if the multiplication $S\otimes_RS\to S$ is an isomorphism and  $\Tor_n^R(S, S)=0$ for all $n\geq 1.$
In this case, we always identify $S\Modcat$ with $\Img(\lambda_*)$, and $\D{S}$ with $\Img(\D{\lambda_*})$.

\subsection{Semi-orthogonal decompositions in triangulated categories\label{sect2.2}}

Now, we recall the definition of semi-orthogonal decompositions in triangulated categories.
Note that semi-orthogonal decompositions are also called hereditary torsion pairs in triangulated categories (see \cite[Section 9.1]{neemanbook}, \cite[Chapter I.2]{BI}).

\begin{Def}\label{SOD}
Let $\mathscr{D}$ be a triangulated category with the shift functor $[1]$. A pair $(\mathscr{X}, \mathscr{Y})$ of full subcategories $\mathscr{X}$ and $\mathscr{Y}$ of $\mathscr{D}$ is called a \emph{semi-orthogonal decomposition} of $\mathscr{D}$ if

$(1)$ $\mathscr{X}$ and $\mathscr{Y}$ are triangulated subcategories of $\mathscr{D}$.

$(2)$ $\Hom_{\mathscr{D}}(X, Y)=0$ for all $X\in\mathscr{X}$ and $Y\in\mathscr{Y}$.

$(3)$ For each object $D\in\mathscr{D}$, there exists a distinguished triangle in $\mathscr{D}$
$$X_D\lra D\lra Y^D\lra X_D[1]$$
with $X_D\in\mathscr{X}$ and $Y^D\in\mathscr{Y}$.

This decomposition is denoted by $\mathscr{D}=<\mathscr{X}, \mathscr{Y}>$.
\end{Def}

Semi-orthogonal decompositions are closely related to half recollements of triangulated categories.
In fact, a pair $(\mathscr{X}, \mathscr{Y})$ of full triangulated subcategories of $\mathscr{D}$ is a semi-orthogonal decomposition of $\mathscr{D}$ if and only if there exists a \emph{lower half recollement} among $\mathscr{X}$, $\mathscr{D}$ and $\mathscr{Y}$, in the sense that
there are four triangle functors demonstrated in the diagram
$$\xymatrix{\mathscr{X}\ar^-{\bf{i}}@/^0.8pc/[r]&\mathscr{D}
\ar^-{\bf{R}}@/^0.8pc/[l]\ar^-{\bf{L}}@/^0.8pc/[r]
&\mathscr{Y}\ar^-{\bf{j}}@/^0.8pc/[l]}$$ such that

$(1)$ ${\bf i}$ and ${\bf j}$ are canonical inclusions;

$(2)$ both $({\bf i}, {\bf R})$ and
$({\bf L}, {\bf j})$ are adjoint pairs;

$(3)$ ${\bf L} {\bf i}=0$ (and thus also
${\bf R}{\bf j}=0$);

$(4)$ for each object $D\in\mathscr{D}$, there exists a distinguished triangle in $\mathscr{D}$
$$
{\bf i}{\bf R}(D)\lra D\lra {\bf j}{\bf L}(D)\lra {\bf i}{\bf R}(D)[1]
$$
where ${\bf i}{\bf R}(D)\ra D$  is the counit adjunction and
$D\ra {\bf j}{\bf L}(D)$ is the unit adjunction. In this case, there are equivalences of triangulated categories: $\mathscr{D}/\mathscr{X}\lraf{\simeq}\mathscr{Y}$ and $\mathscr{D}/\mathscr{Y}\lraf{\simeq}\mathscr{X}.$

Observe that the conditions in Definition \ref{SOD} is weaker than the ones
given in \cite{BGS, Huybrechts, Orlov} because ${\bf i}$ may not have a left adjoint, nor ${\bf j}$  have a right adjoint.
But, if ${\bf i}$ does have a left adjoint (or equivalently, ${\bf L}$ has a fully faithful left adjoint), then the lower half recollement can be completed to a recollement among triangulated categories $\mathscr{X}$, $\mathscr{D}$ and $\mathscr{Y}$ in the sense of Beilinson, Bernstein and Deligne (see \cite{BBD} for definition).

Now, we restate the definition of derived decompositions of abelian categories in a slightly general way.

\begin{Def}\label{fdad}
Let $\mathcal{A}$, $\mathcal{X}$ and $\mathcal{Y}$ be abelian categories, and let $i:\mathcal{X}\to \mathcal{A}$ and $j:\mathcal{Y}\to \mathcal{A}$ be exact functors.
The pair $(\mathcal{X}, \mathcal{Y})$ is called a \emph{derived decomposition} (or $\mathscr{D}^b$-decompostion) of $\mathcal{A}$ with respect to $i$ and $j$ if the two conditions hold:

$(D1)$ The derived functors $\Db{i}: \Db{\mathcal{X}}\to \Db{\mathcal{A}}$ and $\Db{j}: \Db{\mathcal{Y}}\to\Db{\mathcal{A}}$ induced from $i$ and $j$ are fully faithful.

$(D2)'$ $\big(\Img(\Db{i}), \Img(\Db{j})\big)$ is a semi-orthogonal decomposition of $\Db{\mathcal{A}}$.
\end{Def}

Clearly, $(D1)$ implies that both $i$ and $j$ are fully faithful, and thus the images of $i$ and $j$ are abelian subcategories of $\mathcal{A}$.

When $\mathcal{X}$ and $\mathcal{Y}$ are abelian subcategories of $\mathcal{A}$ with inclusions $i:\mathcal{X}\subseteq\mathcal{A}$ and $j:\mathcal{Y}\subseteq\mathcal{A}$, Definition \ref{fdad} and Definition \ref{Intro-def} for $\mathscr{D}^b$-decompositions are equivalent. If $(\mathcal{X}, \mathcal{Y})$ is a derived decomposition of $\mathcal{A}$, then $\mathcal{X}$ and $\mathcal{Y}$ are called \emph{derived factors} of $\mathcal{A}$. Moreover, there are exact sequences of derived categories
$$\xymatrix{\Db{\mathcal{X}}\lraf{\Db{i}}\Db{\mathcal{A}}\lraf{p} \Db{\mathcal{Y}}, \quad \Db{\mathcal{Y}}\lraf{\Db{j}}\Db{\mathcal{A}}\lraf{q} \Db{\mathcal{X}}} $$
where $p$ and $q$ are compositions of the canonical quotient functors with an equivalence functor, respectively.

If $\mathcal{A}$ does not have any non-trivial derived decompositions, then $\mathcal{A}$ is said to be \emph{abelian simple}. Examples of abelian simple categories can eb found in Section \ref{sect4.3}. As in \cite{AKLY} for triangulated categories, we introduce similarly the notion of derived stratifications.

\begin{Def}\label{stratification}
A \emph {derived stratification} of $\mathcal{A}$ is a sequence of
derived decompositions:

$(1)$ a derived decomposition $(\mathcal{X}_0, \mathcal{X}_1)$ of $\mathcal{A}$ if $\mathcal A$
is not abelian simple,

$(2)$ a derived decomposition $(\mathcal{X}_{00}, \mathcal{X}_{01})$ of $\mathcal{X}_0$ if $\mathcal{X}_0$ is not abelian simple,

$(3)$ a derived decomposition $(\mathcal{X}_{10}, \mathcal{X}_{11})$ of $\mathcal{X}_1$ if $\mathcal{X}_1$ is not abelian simple,

$(4)$ derived decompositions of $\mathcal{X}_{ij}$ with $0\le i,j \le 1$ if $\mathcal{X}_{ij}$ are not
abelian simple. Continuing this procedure of decompositions, until one arrives at all derived factors being abelian simple. This procedure may continue to infinitum.
\end{Def}

All the abelian simple categories appearing in this procedure are called \emph {composition factors} of
the stratification. The cardinality of the set of all composition
factors (counting the multiplicity) is called the \emph{length} of the stratification. If this procedure stops after finitely many steps,
we say that this stratification is \emph{finite or of finite length}.

\subsection{Ext-orthogonal pairs and cotorsion pairs in abelian categories}\label{sect2.3}

Derived decompositions of abelian categories are associated with both complete cotorsion pairs and complete Ext-orthogonal pairs in abelian categories.
The notion of complete cotorsion pairs  is classical and has been widely applied to relative homological algebra and generalized tilting theory (see \cite{ EJ, BI, Hovey}), while the notion of complete Ext-orthogonal pairs seems only to be employed in dealing with the telescope conjecture for hereditary rings (see \cite{Krause-S}). We will show in the next section that the latter may be useful in derived decompositions.

Throughout this section, $\mathcal{A}$ is an abelian category, and $(\mathcal{X},\mathcal{Y})$ is a pair of full subcategories of $\mathcal{A}$.

\begin{Def}
$(1)$ The pair $(\mathcal{X}, \mathcal{Y})$ is called a \emph{cotorsion pair} in $\mathcal{A}$ if $\mathcal{X}={^{\bot >0}}\mathcal{Y}$ and $\mathcal{Y}=\mathcal{X}^{\bot >0}$; and a \emph{complete cotorsion pair} in $\mathcal{A}$ if it is a cotorsion pair and, for each object $M\in\mathcal{A}$, there are short exact sequences in $\mathcal{A}$
$$ 0\lra Y_M\lra X_M\lra M\lra 0 \quad \mbox{and}\quad 0\lra M\lra Y^M\lra X^M\lra 0$$
with $X_M, X^M\in\mathcal{X}$ and $Y_M, Y^M\in\mathcal{Y}$.

$(2)$ {\rm \cite[Definition 2.1]{Krause-S} } \label{Ext-orthogonal}
The pair $(\mathcal{X}, \mathcal{Y})$ is said to be \emph{Ext-orthogonal} in $\mathcal{A}$ if $\mathcal{X}={^\bot}\mathcal{Y}$ and $\mathcal{Y}=\mathcal{X}^\bot$; and \emph{complete Ext-orthogonal} in $\mathcal{A}$ if it is Ext-orthogonal and satisfies the gluing condition

\emph{(GC):}\;\; For each object $M\in\mathcal{A}$, there exists a five-term exact sequence in $\mathcal{A}$
$$\varepsilon_M:\;\; 0\lra Y_M\lra X_M\lra M\lra Y^M\lra X^M\lra 0$$
with $X_M, X^M\in\mathcal{X}$ and $Y_M, Y^M\in\mathcal{Y}$.
\end{Def}

\begin{Lem}\label{abelian}
Let $\mathcal{A}$ be an abelian category, and let $\mathcal{X}$ and $\mathcal{Y}$ be full subcategories of $\mathcal{A}$.
If $\mathcal{X}\subseteq {^\bot}\mathcal{Y}$ and $(\mathcal{X}, \mathcal{Y})$ satisfies $(GC)$, then

$(1)$ $\mathcal{X}={^{\bot 0}}\mathcal{Y}\cap {^{\bot 1}}\mathcal{Y}$ and
      $\mathcal{Y}=\mathcal{X}^{\bot 0} \cap \mathcal{X}^{\bot 1}$.

$(2)$ $(\mathcal{X},\mathcal{Y})$ is complete Ext-orthogonal.

$(3)$ Both $\mathcal{X}$ and $\mathcal{Y}$ are abelian subcategories of $\mathcal{A}$.
\end{Lem}

{\it Proof.} $(1)$ It follows from $\mathcal{X}\subseteq {^\bot}\mathcal{Y}$ that  $\mathcal{X}\subseteq {^{\bot 0}}\mathcal{Y}\cap {^{\bot 1}}\mathcal{Y}$. Now, we show the converse of this inclusion. Let $M\in {^{\bot 0}}\mathcal{Y}\cap {^{\bot 1}}\mathcal{Y}$. Since $(\mathcal{X}, \mathcal{Y})$ satisfies $(GC)$, there is a five-term exact sequence $\varepsilon_M$ for $M$ in $\mathcal{A}$. In particular, both $Y_M$ and $Y^M$ belong to $\mathcal{Y}$. It then follows from $\Hom_\mathcal{A}(M, Y^M)=0$ that there is a short exact sequence $0\to Y_M\to X_M\to M\to 0$ in $\mathcal{A}$. Further, this sequence splits due to $\Ext_\mathcal{A}^1(M, Y_M)=0$. Thus $X_M\simeq Y_M\oplus M$. Since  $\Hom_\mathcal{A}(X_M, Y_M)=0$, we have $Y_M=0$ and $M\simeq X_M\in\mathcal{X}$. So ${^{\bot 0}}\mathcal{Y}\cap {^{\bot 1}}\mathcal{Y}\subseteq\mathcal{X}$. Thus the first equality in $(1)$ holds. Similarly, one can verify the second equality in $(1)$.

$(2)$ It suffices to show both $\mathcal{X}={^\bot}\mathcal{Y}$ and $\mathcal{Y}=\mathcal{X}^\bot$. But this follows from $(1)$ and the inclusion $\mathcal{X}\subseteq {^\bot}\mathcal{Y}$.

$(3)$  We only prove that $\mathcal{X}$ is an abelian subcategory of $\mathcal{A}$. The conclusion
on $\mathcal{Y}$ can be proved dually.

Clearly, $\mathcal{X}$ is closed under extensions, kernels of
epimorphisms and cokernels of monomorphisms in $\mathcal{A}$. Thus $\mathcal{X}$ is an abelian subcategory of $\mathcal{A}$ if and only if $\mathcal{X}$ is closed under cokernels (or equivalently, kernels) in
$\mathcal{A}$. In the following, we show that $\mathcal{X}$ is closed under cokernels in $\mathcal{A}$.

Let $f:M\to N$ be a morphism in $\mathcal{A}$ with $M, N\in\mathcal{X}$. Then there is a canonical four-term  exact sequence
$$0\lra \Ker(f)\lra M\lraf{f} N\lra \Coker(f)\lra 0$$
in $\mathcal{A}.$ On the one hand, from $M, N\in\mathcal{X}\subseteq {^{\bot 0}}\mathcal{Y}$  and the fact that ${^{\bot 0}}\mathcal{Y}$ is closed under quotients in $\mathcal{A}$, it follows that
both $\Img(f)$ and $\Coker(f)$ lies in ${^{\bot 0}}\mathcal{Y}$. On the other hand, for $Y\in\mathcal{Y}$, by applying $\Hom_\mathcal{A}(-, Y)$ to the short exact sequence $0\to\Img(f)\to N\to\Coker(f)\to 0$, with $N\in\mathcal{X} ={^\bot}\mathcal{Y}$, we get $\Hom_\mathcal{A}(\Img(f), Y)\simeq \Ext^1_\mathcal{A}(\Coker(f), Y)$. This yields $\Ext^1_\mathcal{A}(\Coker(f), Y)=0$, and consequently, $\Coker(f)\in{^{\bot 0}}\mathcal{Y}\cap {^{\bot 1}}\mathcal{Y}$. By $(1)$,  $\Coker(f)\in\mathcal{X}$ and thus $\mathcal{X}$ is closed under cokernels in $\mathcal{A}$. $\square$

\medskip
Ext-orthogonal pairs have the following properties.

\begin{Lem}{\rm \cite[Lemma 2.9]{Krause-S}}\; \label{elementary}
Let $(\mathcal{X}, \mathcal{Y})$ be an Ext-orthogonal pair in an abelian category $\mathcal{A}$ and $M$ an object in $\mathcal{A}$. Suppose that
there is an exact sequence
$$\varepsilon_M: 0\lra Y_M\lraf{\varepsilon_M^{-2}} X_M\lraf{\varepsilon_M^{-1}} M
\lraf{\varepsilon_M^{0}} Y^M\lraf{\varepsilon_M^{1}} X^M\lra 0$$
in $\mathcal{A}$ with $X_M, X^M\in\mathcal{X}$ and $Y_M, Y^M\in\mathcal{Y}$.

$(1)$ There are isomorphisms of abelian groups for all $X\in\mathcal{X}$ and $Y\in\mathcal{Y}$:
$$
(\varepsilon_M^{-1})^*: \Hom_\mathcal{A}(X, X_M)\lraf{\simeq} \Hom_\mathcal{A}(X, M)\quad\mbox{and}\quad
(\varepsilon_M^{0})_*: \Hom_\mathcal{A}(Y^M, Y)\lraf{\simeq} \Hom_\mathcal{A}(M, Y).
$$

$(2)$ If
$\varepsilon_N:\; 0\ra Y_N\lraf{\varepsilon_N^{-2}} X_N\lraf{\varepsilon_N^{-1}} N
\lraf{\varepsilon_N^{0}} Y^N\lraf{\varepsilon_N^{1}} X^N\ra 0
$
is an exact sequence in $\mathcal{A}$ with $X_N, X^N\in\mathcal{X}$ and $Y_N, Y^N\in\mathcal{Y}$, then each morphism $f: M\to N$ extends uniquely to a morphism $\varepsilon_f: \varepsilon_M\to \varepsilon_N$
of exact sequences:
$$ \xymatrix{
\varepsilon_M: \ar[d]_-{\varepsilon_f} & 0\ar[r] &  Y_M \ar[r]^-{\varepsilon_M^{-2}}\ar[d]_-{Y_f}
& X_M\ar[r]^-{\varepsilon_M^{-1}} \ar[d]_-{X_f}& M\ar[r]^-{\varepsilon_M^{0}}\ar[d]_-{f}
& Y^M\ar[r]^-{\varepsilon_M^{1}} \ar[d]_-{Y^f} & X^M \ar[r] \ar[d]_-{X^f}& 0\\
\varepsilon_N: &0\ar[r] &  Y_N \ar[r]^-{\varepsilon_N^{-2}}
& X_N\ar[r]^-{\varepsilon_N^{-1}} & N \ar[r]^-{\varepsilon_N^{0}}
& Y^N \ar[r]^-{\varepsilon_N^{1}} & X^N \ar[r] & 0. }
$$
$(3)$ Any exact sequence $0\to Y\to X\to M\to Y'\to X'\to 0$ in $\mathcal{A}$ with
$X, X'\in\mathcal{X}$ and $Y, Y'\in\mathcal{Y}$ is isomorphic to $\varepsilon_M$.
\end{Lem}

Now, let $(\mathcal{X},\mathcal{Y})$ be a complete Ext-orthogonal pair in an abelian category $\mathcal{A}$.
For each object $M\in\mathcal{A}$, we fix an exact sequence in $\mathcal{A}$
$$
(\ast)\quad\quad
\varepsilon_M:\quad  0\lra Y_M\lraf{\varepsilon_M^{-2}} X_M\lraf{\varepsilon_M^{-1}} M
\lraf{\varepsilon_M^{0}} Y^M\lraf{\varepsilon_M^{1}} X^M\lra 0
$$
such that $X_M, X^M\in\mathcal{X}$ and $Y_M, Y^M\in\mathcal{Y}$.
In particular, if $M\in\mathcal{X}$, then $\varepsilon_M^{-1}:X_M\to M$ is an isomorphism and $Y^M\simeq 0$; if $M\in\mathcal{Y}$, then $\varepsilon_M^{0}:M\to Y^M$ is an isomorphism and $X_M\simeq 0$.

By Lemma \ref{abelian}, both $\mathcal{X}$ and $\mathcal{Y}$ are abelian subcategories of $\mathcal{A}$ closed under direct summands. Let $i:\mathcal{X}\to \mathcal{A}$ and $j:\mathcal{Y}\to\mathcal{A}$ be the inclusions. Then $i$ and $j$ are exact functors. Moreover, $i$ has a right adjoint $r:\mathcal{A}\to\mathcal{X}$ and $j$ has a left adjoint $\ell:\mathcal{A}\to\mathcal{Y}$, which are defined as follows:

For each $M\in\mathcal{A}$ and for a morphism $f:M\to N$ in $\mathcal{A}$,
$$r(M)=X_M,\;\; r(f)=X_f: X_M\to X_N \quad \mbox{and}\quad \ell(M)=Y^M,\;\; \ell(f)=Y^f: Y^M\to Y^N.$$
This is well defined by Lemma \ref{elementary}. For the adjoint pair $(i,r)$ of functors, the unit adjunction of $X\in\mathcal{X}$ is given by the inverse of the isomorphism $\varepsilon_X^{-1}: r(X)\to X$, and the counit adjunction of $M\in\mathcal{A}$ is given by $\varepsilon_M^{-1}: ir(M)\to M$. Similarly, the unit and counit adjunctions associated with $(\ell, j)$ can be defined by $\varepsilon_M^{0}$.

Now, we can form the following diagram of functors between abelian categories:
$$(\sharp)\quad
\xymatrix{\mathcal{X}\ar^-{i}@/^0.8pc/[r]&\mathcal{A}
\ar^-{r}@/^0.8pc/[l]\ar^-{\ell}@/^0.8pc/[r]
&\mathcal{Y}\ar^-{j}@/^0.8pc/[l],}$$
where $r$ is left exact and $\ell$ is right exact. In general, neither $r$ nor $\ell$ is necessarily exact. So $(\sharp)$ is neither a localization sequence nor a colocalization sequence of abelian categories, and therefore, it may not be completed into a recollement of abelian categories. However, since $i$ and $j$ are exact, they induce derived functors between bounded derived categories:
$\Db{i}:\Db{\mathcal{X}}\to\Db{\mathcal{A}}$ and $\Db{j}:\Db{\mathcal{Y}}\to \Db{\mathcal{A}}$.
With the notation in $(\sharp)$ , the sequence $(\ast)$ can be rewritten as follows:
$$
(\ast)\quad\quad
\varepsilon_M:\quad  0\lra Y_M\lraf{\varepsilon_M^{-2}} r(M)\lraf{\varepsilon_M^{-1}} M
\lraf{\varepsilon_M^{0}} \ell(M)\lraf{\varepsilon_M^{1}} X^M\lra 0.
$$

Finally, we consider the following full subcategories of $\mathcal{A}$ defined via the sequence $(\ast)$:
$$
\mathcal{A}_{\,r\mbox{-}adj}:=\{M\in\mathcal{A}\mid X^M=0\}\quad\mbox{and}\quad \mathcal{A}_{\ell\mbox{-}adj}:=\{M\in\mathcal{A}\mid Y_M=0\}.
$$
The two subcategories have the following properties that will be used in the proof of Theorem \ref{main-result}(1).

\begin{Lem}\label{r-adj}
$(1)$ $\mathcal{A}_{\,r\mbox{-}adj}$ is closed under extensions and quotients in $\mathcal{A}$.

$(2)$ $\mathcal{A}_{\ell\mbox{-}adj}$ is closed under extensions and subobjects in $\mathcal{A}$.

$(3)$ The restriction of $r$ to $\mathcal{A}_{\,r\mbox{-}adj}$ is exact, that is, if $0\to M^{-2}\to M^{-1}\to M^0\to 0$ is an exact sequence in $\mathcal{A}$
with $M^i\in\mathcal{A}_{\,r\mbox{-}adj}$ for $-2\leq i\leq 0$, then $0\to r(M^{-2})\to r(M^{-1})\to r(M^0)\to 0$ is an exact sequence in $\mathcal{X}$.

$(4)$ The restriction of $\ell$ to $\mathcal{A}_{\ell\mbox{-}adj}$ is exact, that is, if $0\to M^{-2}\to M^{-1}\to M^0\to 0$ is an exact sequence in $\mathcal{A}$
with $M^i\in\mathcal{A}_{\ell\mbox{-}adj}$ for $-2\leq i\leq 0$, then $0\to \ell(M^{-2})\to \ell(M^{-1})\to \ell(M^0)\to 0$ is an exact sequence in $\mathcal{Y}$.

\end{Lem}

{\it Proof.}  We only  prove $(1)$ and $(3)$ since $(2)$ and $(4)$ can be proved dually.

Let $0\to M^{-2}\lraf{f} M^{-1}\lraf{g} M^0\to 0$ be an exact sequence in $\mathcal{A}$. We regard it as a complex $\cpx{M}$ in $\Cb{\mathcal{A}}$ with $M^0$ in degree $0$. It follows from Lemma \ref{elementary}(2) that the sequence $(\ast)$, associated with $M^i$ for $-2\leq i\leq 0$, induces an exact sequence of complexes over $\mathcal{A}$:
$$
0\lra Y_{\cpx{M}}\lraf{\varepsilon_{\cpx{M}}^{-2}}r(\cpx{M})\lraf{\varepsilon_{\cpx{M}}^{-1}} \cpx{M}\lraf{\varepsilon_{\cpx{M}}^{0}}\ell(\cpx{M})
\lraf{\varepsilon_{\cpx{M}}^{1}} X^{\cpx{M}}\to 0,
$$
where $r(\cpx{M})$ means $(r(M^i))_{i\in \mathbb{Z}}$. Recall that $r$ is a left exact functor and $\ell$ is a right exact functor. Thus the complexes $r(\cpx{M})$ and $\ell(\cpx{M})$ are exact everywhere except in the degrees $0$ and $-2$,
respectively. Since $\cpx{M}$ is an exact sequence, $(1)$ holds.

(3) Suppose $M^i\in\mathcal{A}_{\,r\mbox{-}adj}$ for $-2\leq i\leq 0$. Then $X^{\cpx{M}}=0$. To show that $r(\cpx{M})$ is a exact sequence, it suffices to show that the homomorphism $r(g): r(M^{-1})\to r(M^0)$ is surjective
(or equivalently, $H^0(r(\cpx{M}))=0$).

Let $c(\cpx{M})$ be the cokernel of the chain map $\varepsilon_{\cpx{M}}^{-2}$. Taking cohomologies on the sequence $0\to Y_{\cpx{M}}\to r(\cpx{M})\to c(\cpx{M})\to 0$, we get a long exact sequence in $\mathcal{A}$:
$$H^{-1}(r(\cpx{M}))\lra H^{-1}(c(\cpx{M}))\lra H^0(Y_{\cpx{M}})\lra H^0(r(\cpx{M}))\lra H^0(c(\cpx{M})).$$
Note that $H^{-1}(r(\cpx{M}))=0$ and $H^0(Y_{\cpx{M}})\simeq \Coker(Y_g)$, where $Y_g: Y_{M^{-1}}\to Y_{M^0}$ is induced from $g$.
Since the sequence $\cpx{M}$ is exact and $0\to c(\cpx{M})\to\cpx{M}\to\ell(\cpx{M})\to 0$ is exact
in $\Cb{\mathcal{A}}$, we have $H^i(c(\cpx{M}))\simeq H^{i-1}(\ell(\cpx{M}))$ for any $i\in\mathbb{Z}$. This implies
$H^0(c(\cpx{M}))\simeq H^{-1}(\ell(\cpx{M}))=0$ and $H^{-1}(c(\cpx{M}))\simeq H^{-2}(\ell(\cpx{M}))\simeq \Ker(Y^f)$, where
$Y^f: \ell(M^{-2})\to \ell(M^{-1})$ is induced from $f$. Consequently, there is a short exact sequence in $\mathcal{A}$:
$$ 0\lra \Ker(Y^f)\lra \Coker(Y_g)\lra H^0(r(\cpx{M}))\lra 0.$$
Clearly, $\Ker(Y^f), \Coker(Y_g)\in\mathcal{Y}$ and $H^0(r(\cpx{M}))\in\mathcal{X}$ since $\mathcal{Y}$ and $\mathcal{X}$ are abelian full subcategories of $\mathcal{A}$ by Lemma \ref{abelian}(3). This implies $H^0(r(\cpx{M}))\in\mathcal{X}\cap\mathcal{Y}$. It follows from $\mathcal{X}={^\bot}\mathcal{Y}$ that $H^0(r(\cpx{M}))=0$. Thus $(3)$ holds. $\square$

\section{Derived decompositions of abelian categories\label{sect3}}
In this section we will prove Theorem \ref{main-result}. In particular, we show that a complete Ext-orthogonal pair is a derived decomposition of an abelian category
with enough projectives and injectives if and only if the five-term exact sequences for both projective and injective objects are reduced to
four terms.



\subsection{Proof of Theorem \ref{main-result}(1)}

Complete Ext-orthogonal pairs and derived decompositions, both are defined at the level of abelian categories. But the latter
reflect information on bounded derived categories of abelian categories. This suggests that derived decompositions might imply complete Ext-orthogonal pairs. In the following we will show this implication.

\begin{Prop}\label{derived}
Let $\mathcal{A}$ be an abelian category and let $\mathcal{X}$ and $\mathcal{Y}$ be full subcategories of $\mathcal{A}$. Given $*\in \{b, +,-, \varnothing\}$, if $(\mathcal{X}, \mathcal{Y})$ is a $\mathscr{D}^*$-decomposition of $\mathcal{A}$, then it is a complete Ext-orthogonal pair in $\mathcal{A}$.
\end{Prop}

{\it Proof.} Suppose that $(\mathcal{X}, \mathcal{Y})$ is a derived decomposition of $\mathcal{A}$ in $\Ds{\mathcal{A}}$. By Definition \ref{Intro-def}(D2), we have $\mathcal{X}\subseteq {^\bot}\mathcal{Y}$. Moreover, by Definition \ref{Intro-def}(D3), for each $M\in\mathcal{A}$, there is a triangle $\cpx{X}\to M\to \cpx{Y}\to \cpx{X}[1]$ in $\Ds{\mathcal{A}}$ such that $\cpx{X}\in\Ds{\mathcal{X}}$ and $\cpx{Y}\in\Ds{\mathcal{Y}}$, where $\Ds{\mathcal{X}}$ and $\Ds{\mathcal{Y}}$ can be regarded as triangulated subcategories of $\Db{\mathcal{A}}$ by $(D1)$ in Definition \ref{Intro-def}.
By taking cohomologies on this triangle, one gets the following long exact sequence in $\mathcal{A}$:
$$
0\to H^{-1}(\cpx{Y})\to H^0(\cpx{X})\to M\to H^0(\cpx{Y})\to H^1(\cpx{X})\to 0.
$$
Recall that $\mathcal{X}$ and $\mathcal{Y}$ are abelian subcategories of $\mathcal{A}$.
Since $\cpx{X}\in\Ds{\mathcal{X}}$, we have $H^i(\cpx{X})\in\mathcal{X}$ for all $i\in\mathbb{Z}$. In particular, both $H^0(\cpx{X})$ and $H^1(\cpx{X})$ lie in $\mathcal{X}$. Similarly, one can show that both $H^{-1}(\cpx{Y})$ and $H^0(\cpx{Y})$ belong to $\mathcal{Y}$. Thus the pair $(\mathcal{X},\mathcal{Y})$ satisfies the gluing condition. Now, it follows from Lemma \ref{abelian}(2) that $(\mathcal{X},\mathcal{Y})$ is a complete Ext-orthogonal pair in $\mathcal{A}$. $\square$

\medskip
Having shown that $\mathscr{D}^*$-decompositions are complete Ext-orthogonal pairs, we now consider the converse question:

\medskip
\emph{ Given a complete Ext-orthogonal pair $(\mathcal{X}, \mathcal{Y})$ in an abelian category $\mathcal{A}$, when is it a $\mathscr{D}^*$-decomposition of $\mathcal{A}$ for $*\in \{b, +,-, \varnothing\}?$ }

\medskip
To address this question, we use derived categories of exact categories.
Recall that an exact category (in the sense of Quillen) is an additive category $\mathcal{E}$ endowed with a class of conflations
closed under isomorphism and satisfying certain axioms (for example, see \cite[Section 4]{Keller}). In case that $\mathcal{E}$ is an abelian
category, the class of conflations coincides with the class of short exact sequences.

Let $\mathcal{E}$ be a full subcategory of the abelian category $\mathcal{A}$.
Suppose that $\mathcal{E}$ is closed under extensions in $\mathcal{A}$,
that is, for any exact sequence $0\to X\to Y\to Z\to 0$ in $\mathcal{A}$ with both $X,Z\in \mathcal{E}$, we have $Y\in\mathcal{E}$. Then $\mathcal{E}$ endowed with the short exact sequences of $\mathcal{A}$ having their terms in $\mathcal{E}$ is an exact
category and the inclusion $\mathcal{E}\subseteq\mathcal{A}$ is a fully faithful exact functor. So, $\mathcal{E}$ is called a
\emph{fully exact subcategory} of $\mathcal{A}$.

A complex $\cpx{X}\in\C{\mathcal{E}}$ is said to be \emph{strictly exact} if it is exact in $\C{\mathcal{A}}$ and all of its boundaries belong to $\mathcal{E}$.
Let $\mathscr{K}_{\rm {ac}}(\mathcal{E})$ be the full subcategory of $\K{\mathcal{E}}$ consisting of those complexes which are isomorphic to
strictly exact complexes. Then $\mathscr{K}_{\rm {ac}}(\mathcal{E})$ is a full triangulated subcategory of $\K{\mathcal{E}}$.
The \emph{unbounded derived category} of $\mathcal{E}$, denoted by $\D{\mathcal{E}}$, is defined to be the Verdier quotient of $\K{\mathcal{E}}$
by $\mathscr{K}_{\rm {ac}}(\mathcal{E})$. Similarly, the bounded-below, bounded-above and bounded derived categories
$\Dz{\mathcal{E}}$, $\Df{\mathcal{E}}$ and $\Db{\mathcal{E}}$ can be defined through bounded-below, bounded-above and bounded
complexes over $\mathcal{E}$, respectively. Moreover, the canonical functor $\mathscr{D}^*(\mathcal{E})\to\D{\mathcal{E}}$ is fully faithful
for any $*\in\{+, -, b\}$. Note that if $\mathcal{E}$ is closed under cokernels of monomorphisms in $\mathcal{A}$, then
$\cpx{X}\in\Cb{\mathcal{E}}$ is strictly exact if and only if it is exact in $\C{\mathcal{A}}$.

For more details on derived categories of exact categories, we refer the reader to \cite{Keller}.
The following result follows from \cite[Theorem 12.1]{Keller}.

\begin{Lem}\label{resolution}
Let $\mathcal{E}$ be a full subcategory of an abelian category $\mathcal{A}$. Assume that the two conditions hold:

$(a)$ If $0\to X\to Y\to Z\to 0$ is an exact sequence in $\mathcal{A}$ with $X\in \mathcal{E}$, then $Y\in\mathcal{E}$ if and only if $Z\in\mathcal{E}$.

$(b)$ For each object $M\in\mathcal{A}$, there is a long exact sequence in $\mathcal{A}$:
$$
0\lra M\lra  E_0\lra E_{1}\lra \cdots\lra E_{n-1} \lra E_n\lra 0
$$
\quad\quad for a nonnegative integer $n$ such that $E_i\in\mathcal{E}$ for all $0\leq i\leq n$.

\smallskip
\noindent Then the inclusion $\mathcal{E}\subseteq\mathcal{A}$ induces a triangle equivalence $\Dz{\mathcal{E}}\lraf{\simeq} \Dz{\mathcal{A}}$
which can be restricted to an equivalence $\Db{\mathcal{E}}\lraf{\simeq} \Db{\mathcal{A}}$.
\end{Lem}

\begin{Rem}\label{unbounded resolution}
In Lemma \ref{resolution}, the inclusion $\mathcal{E}\subseteq\mathcal{A}$ induces a triangle equivalence $\mathscr{D}(\mathcal{E})\lraf{\simeq} \mathscr{D}(\mathcal{A})$ which can be restricted to an equivalence $\mathscr{D}^*(\mathcal{E})\lraf{\simeq} \mathscr{D}^*(\mathcal{A})$ for any
$*\in\{+, -, b\}$ provided that $(b)$ in Lemma \ref{resolution} is strengthened by

($b'\,$) There is a nonnegative integer $n$ such that, for each object $M\in\mathcal{A}$, there is a long exact sequence
$0\to M\to E_0\to E_{1}\to \cdots\to E_{n-1}\to E_n\to 0$
in $\mathcal{A}$ with $E_i\in\mathcal{E}$ for all $0\leq i\leq n$.

\medskip
For a proof of this result and its dual statement, we refer the reader to \cite[Proposition A.5.6]{Pos1}
\end{Rem}

\noindent {\bf Proof of Theorem \ref{main-result}(1).}
By $(a)$ and $(b)$, it follows from Lemma \ref{abelian} that $(\mathcal{X},\mathcal{Y})$ is a complete Ext-orthogonal pair in $\mathcal{A}$
and that both $\mathcal{X}$ and $\mathcal{Y}$ are abelian subcategories of $\mathcal{A}$. In particular, $(D2)$ in Definition \ref{Intro-def} is satisfied.

Now, we keep all the notation introduced in Section \ref{sect2.3}.
Under the assumptions of $(a)$, $(b)$ and $(c)$, we show that the functor $\Ds{i}:\Ds{\mathcal{X}}\to\Ds{\mathcal{A}}$, induced from the inclusion $i:\mathcal{X}\to\mathcal{A}$, is fully faithful.

By Lemma \ref{r-adj}(1), $\mathcal{A}_{\,r\mbox{-}adj}$ is closed under extensions in $\mathcal{A}$ and thus a fully exact subcategory of $\mathcal{A}$. Moreover, $i$ has a right adjoint $r:\mathcal{A}\to\mathcal{X}$ which is an exact functor when restricted to $\mathcal{A}_{\,r\mbox{-}adj}$ by Lemma \ref{r-adj}(3).
Thus $r:\mathcal{A}_{\,r\mbox{-}adj}\to \mathcal{X}$ induces a derived functor $\Ds{r}: \Ds{\mathcal{A}_{\,r\mbox{-}adj}}\to \Ds{\mathcal{X}}$.
Since the functor $i$ is fully faithful, the composition of $i$ and $r$ is isomorphic to the identity functor of $\mathcal{X}$. Thus $(\Ds{i}, \Ds{r})$ is
an adjoint pair and the composition of $\Ds{i}$ with $\Ds{r}$ is isomorphic to the identity functor of $\Ds{\mathcal{X}}$. This implies that $\Ds{i}:\Ds{\mathcal{X}}\to \Ds{\mathcal{A}_{\,r\mbox{-}adj}}$ is fully faithful. Further, for each object $M\in\mathcal{A}$, it follows from $(c)$ that there is a monomorphism $M\to I$ in $\mathcal{A}$ such that $I\in \mathcal{A}_{\,r\mbox{-}adj}$. Since $\mathcal{A}_{\,r\mbox{-}adj}$ is closed under quotients in $\mathcal{A}$ by Lemma \ref{r-adj}(1), there is an exact sequence $0\to M\to I\to J\to 0$ in $\mathcal{A}$ with $I, J\in\mathcal{A}_{\,r\mbox{-}adj}$. By Remark \ref{unbounded resolution}, the inclusion $\mathcal{A}_{\,r\mbox{-}adj}\subseteq \mathcal{A}$ of exact categories induces a triangle equivalence $\mathscr{D}^*(\mathcal{A}_{\,r\mbox{-}adj})\lraf{\simeq} \mathscr{D}^*(\mathcal{A})$ for any $*\in \{b, +,-, \varnothing\}.$
Consequently, $\Ds{i}:\Ds{\mathcal{X}}\to \Ds{\mathcal{A}}$ is fully faithful.

Dually, under the assumptions of $(a)$, $(b)$ and $(d)$, the functor $\Ds{j}:\Ds{\mathcal{Y}}\to\Ds{\mathcal{A}}$, induced from the inclusion $j:\mathcal{Y}\to\mathcal{A}$, is fully faithful. This follows from Lemma \ref{r-adj}(2), Lemma \ref{r-adj}(4) and the dual statement of Lemma \ref{resolution}.

Thus Definition \ref{Intro-def}(D1) is satisfied. Now, we identify $\Ds{\mathcal{X}}$ and $\Ds{\mathcal{Y}}$ with $\Img(\Ds{i})$ and $\Img(\Ds{j})$, respectively. It remains to check Definition \ref{Intro-def}(D3).

Let $\cpx{N}\in\Ds{\mathcal{A}}$. Since $\Ds{\mathcal{A}_{\,r\mbox{-}adj}}\lraf{\simeq} \Ds{\mathcal{A}}$, there is a complex $\cpx{M}:=(M^i)_{i\in\mathbb{Z}}\in \mathscr{C}^*(\mathcal{A}_{\,r\mbox{-}adj})$ such that $\cpx{N}\simeq\cpx{M}$ in $\Ds{\mathcal{A}}$. Moreover, since $(\mathcal{X},\mathcal{Y})$ is a complete Ext-orthogonal pair in $\mathcal{A}$, it follows from Lemma \ref{elementary}(2) that each morphism $f: M\to N$ in $\mathcal{A}$ extends uniquely to a morphism $\varepsilon_f: \varepsilon_M\to \varepsilon_N$ of exact sequences (see Lemma \ref{elementary} for notation). Applying this to the differentials of $\cpx{M}$ yields a long exact sequence
$$ \xymatrix{
\varepsilon_{\cpx{M}}:& 0\ar[r] &  Y_{\cpx{M}} \ar[r]^-{\varepsilon_{\cpx{M}}^{-2}}
& X_{\cpx{M}}\ar[r]^-{\varepsilon_{\cpx{M}}^{-1}} & \cpx{M}\ar[r]^-{\varepsilon_{\cpx{M}}^{0}}
& Y^{\cpx{M}}\ar[r]^-{\varepsilon_{\cpx{M}}^{1}}& X^{\cpx{M}} \ar[r]& 0}
$$
in $\mathscr{C}^*(\mathcal{A})$ such that $X_{\cpx{M}}, X^{\cpx{M}}\in\mathscr{C}^*(\mathcal{X})$ and $Y_{\cpx{M}}, Y^{\cpx{M}}\in\mathscr{C}^*(\mathcal{Y})$.
Since $M^i\in \mathcal{A}_{\,r\mbox{-}adj}$ for any $i\in\mathbb{Z}$, we have
$X^{M^i}=0$. This implies $X^{\cpx{M}}=0$ and thus the exact sequence $\varepsilon_{\cpx{M}}$ is of the form:
$$ \xymatrix{
\varepsilon_{\cpx{M}}:& 0\ar[r] &  Y_{\cpx{M}} \ar[r]^-{\varepsilon_{\cpx{M}}^{-2}}
& X_{\cpx{M}}\ar[r]^-{\varepsilon_{\cpx{M}}^{-1}} & \cpx{M}\ar[r]^-{\varepsilon_{\cpx{M}}^{0}}
& Y^{\cpx{M}}\ar[r]& 0.}
$$
Let $Z_{\cpx{M}}$ be the mapping cone of the chain map ${\varepsilon_{\cpx{M}}^{-1}}$. Then $X_{\cpx{M}}\to \cpx{M}\to Z_{\cpx{M}}\to X_{\cpx{M}}[1]$
is a distinguished triangle in $\Ds{\mathcal{A}}$. We show $Z_{\cpx{M}}\in\Ds{\mathcal{Y}}$.

Since $\varepsilon_{\cpx{M}}^{-1}$ is the composite of the surjection $ X_{\cpx{M}}\to \Coker(\varepsilon_{\cpx{M}}^{-2})$ with the injection $\Coker(\varepsilon_{\cpx{M}}^{-2})\to \cpx{M}$, there is a triangle $Y_{\cpx{M}}[1]\to Z_{\cpx{M}}\to Y^{\cpx{M}}\to Y_{\cpx{M}}[2]$ in $\Ds{\mathcal{A}}$ (constructed from the octahedral axiom of triangulated categories).
As $\Ds{\mathcal{Y}}$ is a full triangulated subcategory of $\Ds{\mathcal{A}}$,  it follows from $Y_{\cpx{M}},  Y^{\cpx{M}}\in\Ds{\mathcal{Y}}$ that $Z_{\cpx{M}}\in\Ds{\mathcal{Y}}$.

Since $\cpx{N}\simeq\cpx{M}$ in $\Ds{\mathcal{A}}$, there is a triangle
$X_{\cpx{M}}\to \cpx{N}\to Z_{\cpx{M}}\to X_{\cpx{M}}[1]$ satisfying that $X_{\cpx{M}}\in\Ds{\mathcal{X}}$ and $Z_{\cpx{M}}\in\Ds{\mathcal{Y}}$. This shows Definition \ref{Intro-def}(D3). Thus $(\mathcal{X},\mathcal{Y})$ is a derived decomposition of $\Ds{\mathcal{A}}$. $\square$

\medskip
As a consequence of Theorem \ref{main-result}(1), we have

\begin{Koro}
Let $\mathcal{A}$ be an abelian category, and let $\mathcal{X}$ and $\mathcal{Y}$ be full subcategories of $\mathcal{A}$.
Suppose that $\mathcal{X}\subseteq {^\bot}\mathcal{Y}$ and for each object $M\in\mathcal{A}$, there is a short exact sequence
$0\to X\to M\to Y\to 0$ in $\mathcal{A}$ such that $X\in\mathcal{X}$ and $Y\in\mathcal{Y}$.
Then $(\mathcal{X}, \mathcal{Y})$ is a $\mathscr{D}^*$-decomposition of $\mathcal{A}$ for any
$\ast\in \{b,+,-, \varnothing\}$.
\end{Koro}

\subsection{Proof of Theorem \ref{main-result}(2)}

Throughout this section, $\mathcal{A}$ is an abelian category and $(\mathcal{X},\mathcal{Y})$ is a complete Ext-orthogonal pair in $\mathcal{A}$.
We keep all the notation introduced in Section \ref{sect2.3}.

\begin{Lem}\label{longexact}
Let $M\in\mathcal{A}$ and $N\in\mathcal{Y}$.

$(1)$ There is a long exact sequence of extension groups:
$$
\cdots\lra \Ext_\mathcal{A}^{n-2}(Y_M, N)\lra \Ext_\mathcal{A}^n(\ell{(M)}, N)\lra
\Ext_\mathcal{A}^n(M, N)\lra \Ext_\mathcal{A}^{n-1}(Y_M, N)\lra\cdots
$$
where $n\in\mathbb{Z}$.

$(2)$ If $M\in {^{\bot {1}}}\mathcal{Y}$, then $\ell(M)\in {^{\bot {1}}}\mathcal{Y}$.

$(3)$ If $M\in {^{\bot {>0}}}\mathcal{Y}$, then $\Ext_\mathcal{A}^{n-2}(Y_M, N)\simeq  \Ext_\mathcal{A}^n(\ell(M), N)$ for all $n\geq 2$.
\end{Lem}

{\it Proof.} The sequence in $(1)$ follows from applying $\Ext_\mathcal{A}^n(-, N)$ for $N\in\mathcal{Y}$ to $\varepsilon_M$ and the fact that $r(M)=X_M, X^M\in\mathcal{X}$ and $\mathcal{X}={^\bot}\mathcal{Y}$, while $(2)$ and $(3)$ follow from $(1)$. $\square$

\medskip
From Lemma \ref{longexact}, we have the following

\begin{Koro}\label{projective}
Assume that $\mathcal{A}$ has enough projectives. Let $P\in\mathscr{P}(\mathcal{A})$ and $N\in\mathcal{Y}$.

$(1)$ $\ell(P) \in {^{\bot {1}}}\mathcal{Y}$ and $\Ext_\mathcal{A}^{n-2}(Y_P, N)\simeq  \Ext_\mathcal{A}^n(\ell(P), N)$ for all $n\geq 2$.

$(2)$ $\ell(P)\in {^{\bot {>0}}}\mathcal{Y}$ if and only if $Y_P=0$.
\end{Koro}
{\it Proof.} If $P\in\mathscr{P}(\mathcal{A})$, then $P\in {^{\bot {>0}}}\mathcal{Y}$. Now, $(1)$ follows from Lemma \ref{longexact}(2)-(3). Further, by $(1)$,  $\ell(P)\in {^{\bot {>0}}}\mathcal{Y}$ if and only if $\Ext_\mathcal{A}^{m}(Y_P, Z)=0$ for all $Z\in\mathcal{Y}$ and
for all $m\geq 0$. Due to $Y_P\in\mathcal{Y}$,  $\Hom_\mathcal{A}(Y_P, Y_P)=0$ implies $Y_P=0$. Thus $(2)$ holds. $\square$

\medskip
We need the following result from \cite[Proposition 4.3]{GL}.

\begin{Lem}\label{embedding}
Let $\mathcal{B}$ be an abelian full subcategory of $\mathcal{A}$ and let $\lambda:\mathcal{B}\to\mathcal{A}$ be the inclusion. Then $\Db{\lambda}:\Db{\mathcal{B}}\to\Db{\mathcal{A}}$ is fully faithful if and only if, for any $X, Y\in\mathcal{B}$ and for any $n\in\mathbb{N}$, the homomorphism $\varphi_{X,Y}^n:\Ext_\mathcal{B}^n(X,Y)\to \Ext_\mathcal{A}^n(X,Y)$ induced from $\Db{\lambda}$ is an isomorphism.
\end{Lem}

The following result, which will be used in Section \ref{sect4}, is implied by Lemmas \ref{longexact}(1) and  \ref{embedding}. Here, we omit its proof.

\begin{Koro}\label{proj or injective}
$(1)$ Suppose that $\Db{j}:\Db{\mathcal{Y}}\to \Db{\mathcal{A}}$ is fully faithful. If
$I\in\mathscr{I}(\mathcal{Y})$, then $\Ext_\mathcal{A}^n(M, I)=0$ for any $M\in\mathcal{A}$ and $n\geq 2$.

$(2)$ Suppose that $\Db{i}:\Db{\mathcal{X}}\to \Db{\mathcal{A}}$ is fully faithful. If
$P\in\mathscr{P}(\mathcal{X})$, then $\Ext_\mathcal{A}^n(P, M)=0$ for any $M\in\mathcal{A}$ and $n\geq 2$.
\end{Koro}

\begin{Lem}\label{ff}
Let $\mathcal{B}$ be an abelian full subcategory of $\mathcal{A}$.

$(1)$ Assume that $\mathcal{B}$ has enough projectives. Then the derived functor $\Db{\mathcal{B}}\to\Db{\mathcal{A}}$ induced from the inclusion
$\mathcal{B}\subseteq \mathcal{A}$ is fully faithful if and only if
$\mathscr{P}(\mathcal{B})=\mathcal{B}\cap {^{\bot {>0}}}\mathcal{B}$.

$(2)$ Assume that $\mathcal{B}$ has enough injectives. Then the derived functor $\Db{\mathcal{B}}\to\Db{\mathcal{A}}$ induced from the inclusion
$\mathcal{B}\subseteq \mathcal{A}$ is fully faithful if and only if
$\mathscr{I}(\mathcal{B})=\mathcal{B}\cap\mathcal{B}^{\bot {>0}}$.
\end{Lem}

{\it Proof.} We only  prove $(1)$ since $(2)$ can be proved dually.
Note that the inclusions $\mathcal{B} \cap {^{\bot {>0}}}\mathcal{B}\subseteq \mathcal{B}\cap{^{\bot 1}}\mathcal{B}\subseteq \mathscr{P}(\mathcal{B})$ always hold.

Let $\lambda:\mathcal{B}\to\mathcal{A}$ be the inclusion. By Lemma \ref{embedding}, to prove (1), it is enough to show that, for any $X, Y\in\mathcal{B}$ and for any $n\in\mathbb{N}$, the homomorphism $\varphi_{X,Y}^n:\Ext_\mathcal{B}^n(X,Y)\to \Ext_\mathcal{A}^n(X,Y)$ is an isomorphism. Clearly, $\varphi_{X,Y}^0$ is the identity map since $\mathcal{B}$ is a full subcategory of $\mathcal{A}$. So it suffices to check that  $\varphi_{X,Y}^n$ is an isomorphism for $n\geq 1$.

Suppose that $\Db{\lambda}$ is fully faithful. If $X\in\mathscr{P}(\mathcal{B})$, then $\Ext_\mathcal{B}^n(X,Y)=0$ for all $n\geq 1$, and therefore $\Ext_\mathcal{A}^n(X,Y)=0$.
Thus $X\in {^{\bot {>0}}}\mathcal{B}$ and $\mathscr{P}(B)\subseteq \mathcal{B}\cap {^{\bot {>0}}}\mathcal{B}$. So the necessity of Lemma \ref{ff}(1) holds.

Conversely, suppose $\mathscr{P}(B)=\mathcal{B}\cap {^{\bot {>0}}}\mathcal{B}$. This implies that $\varphi_{X,Y}^n$ is an isomorphism for $X\in\mathscr{P}(\mathcal{B})$ because both $\Ext_\mathcal{B}^n(X,Y)$ and $\Ext_\mathcal{A}^n(X,Y)$ vanish. Since $\mathcal{B}$ has enough projectives, any object  $X\in\mathcal{B}$ has a projective resolution in $\mathcal{B}$. This resolution is also exact in $\mathcal{A}$ because $\mathcal{B}$ is an abelian full subcategory of $\mathcal{A}$. Now, for a fixed object $Y\in\mathcal{B}$, we apply the functors $\Ext_\mathcal{B}^i(-, Y)$  and $\Ext_\mathcal{A}^i(-, Y)$ for $i\in\mathbb{N}$ to this resolution and then get two long exact sequences of extension groups, linked by commutative diagrams. Note that $\Ext_\mathcal{A}^i(P, Y)=0$ for all $i\geq 1$ and $P\in\mathscr{P}(\mathcal{B})$. By induction on $n$ and by the Five-Lemma, we can show that $\varphi_{X,Y}^n$ are isomorphisms for $n\geq 1$. Thus $\Db{\lambda}$ is fully faithful. $\square$

\begin{Lem}\label{transform}
$(1)$ If $\mathcal{A}$ has enough projectives, then so does $\mathcal{Y}$, and
$\mathscr{P}(\mathcal{Y})=\add\big(\{\ell(P)\mid P\in\mathscr{P}(\mathcal{A})\}\big).$

$(2)$ If $\mathcal{A}$ has enough injectives, then so does $\mathcal{X}$, and
$\mathscr{I}(\mathcal{X})=\add\big(\{r(I)\mid I\in\mathscr{I}(\mathcal{A})\}\big).$
\end{Lem}

{\it Proof.} $(1)$
Since $j$ is an exact functor and $(\ell,j)$ is an adjoint pair, $\ell$ is a right exact functor and
preserves projective objects. This means $\ell(P)\in\mathscr{P}(\mathcal{Y})$ for  $P\in\mathscr{P}(\mathcal{A})$. Given any object $Y\in\mathcal{Y}$, since $\mathcal{A}$ has enough projectives, there exists an epimorphism $\pi: Q\to j(Y)$ in $\mathcal{A}$ with $Q\in\mathscr{P}(\mathcal{A})$. Hence $\ell(\pi):\ell(Q)\to \ell(j(Y))$ is an epimorphism in $\mathcal{Y}$. As $\ell(j(Y))\simeq Y$, $\ell(\pi)$ is an epimorphism from $\ell(Q)$ to $Y$. This shows that $\mathcal{Y}$ has enough projectives. Moreover, if $Y\in\mathscr{P}(\mathcal{Y})$, then
$Y$ is a direct summand of $\ell(Q)$. This shows $(1)$.

 $(2)$ can be proved dually. $\square$

\medskip
The next result characterizes when $\Db{i}$ and $\Db{j}$ are fully faithful.

\begin{Lem} \label{EOff}
$(1)$ If $\mathcal{A}$ has enough projectives, then
$\Db{j}:\Db{\mathcal{Y}}\to \Db{\mathcal{A}}$ is fully faithful if and only if $Y_P=0$ for any object $P\in\mathscr{P}(\mathcal{A})$.

$(2)$ If $\mathcal{A}$ has enough injectives, then
$\Db{i}:\Db{\mathcal{X}}\to \Db{\mathcal{A}}$ is fully faithful if and only if $X^I=0$ for any
object $I\in\mathscr{I}(\mathcal{A})$.
\end{Lem}

{\it Proof.} We show $(1)$ by Lemma \ref{ff}(1).
By Lemmas \ref{transform}(1) and \ref{ff}(1), the functor $\Db{j}$ is fully faithful if and only if $\ell(P)\in {^{\bot {>0}}}\mathcal{Y}$ for all $P\in\mathscr{P}(\mathcal{A})$. But the latter is  equivalent to saying $Y_P=0$ by Corollary \ref{projective}(2). Thus $(1)$ holds. Dually, $(2)$ can be proved by Lemma \ref{ff}(2). $\square$

\medskip
\noindent {\bf Proof of Theorem \ref{main-result}(2).}
The sufficiency is a direct consequence of Theorem \ref{main-result}(1). To show the necessity, let $(\mathcal{X}, \mathcal{Y})$ be a $\mathscr{D}^*$-decomposition of $\mathcal{A}$. Then $(a)$ and $(b)$ hold by Proposition \ref{derived}. If $\mathcal{A}$ has enough injectives, then $(c)$ and $c'$ are equivalent. Similarly, if $\mathcal{A}$ has enough projectives, then $(d)$ and $(d')$ are equivalent. Note that $\Ds{i}:\Ds{\mathcal{X}}\to\Ds{\mathcal{A}}$ can be restricted to bounded derived categories. By Definition \ref{Intro-def}(D1), $\Db{i}:\Db{\mathcal{X}}\to\Db{\mathcal{A}}$ is fully faithful. Similarly, $\Db{j}: \Db{\mathcal{Y}}\to \Db{\mathcal{A}}$ is fully faithful. Since $\mathcal{A}$ has enough projectives and injectives, it follows from Lemma \ref{EOff} that $(c')$ and $(d')$ in Theorem \ref{main-result}(2) hold. $\square$

\subsection{$\mathscr{D}^*$-decompositions of $\mathcal{A}$ and semi-orthogonal decompositions of $\Ds{\mathcal{A}}$\label{Section 3.3}}

In this section we establish relations between derived decompositions of abelian categories and semi-orthogonal decompositions of different derived categories.

The following result is an easy observation from Definitions \ref{SOD} and \ref{Intro-def}.
\begin{Lem}\label{DD-HTP}
Let $\mathcal{A}$ be an abelian category, $\mathcal{X}$ and $\mathcal{Y}$ abelian subcategories of $\mathcal{A}$ and $\ast\in \{b,+,-, \varnothing\}$.
Suppose that the inclusions $i: \mathcal{X}\subseteq\mathcal{A}$ and $j: \mathcal{Y}\subseteq\mathcal{A}$ induce fully faithful functors $\Ds{i}:\Ds{\mathcal{X}}\to \Ds{\mathcal{A}}$ and $\Ds{j}:\Ds{\mathcal{Y}}\to\Ds{\mathcal{A}}$, respectively. If $\big(\Img(\Ds{i}), \Img(\Ds{j})\big)$ is a semi-orthogonal decomposition of $\Ds{\mathcal{A}}$, then $(\mathcal{X}, \mathcal{Y})$ is a $\mathscr{D}^*$-decomposition of $\mathcal{A}$.
\end{Lem}

To obtain the converse of Lemma \ref{DD-HTP}, we consider abelian categories with additional properties.
\begin{Def}
An abelian category $\mathcal{A}$

$(1)$ is complete (respectively, cocomplete) if products (respectively, coproducts) indexed over sets exist in $\mathcal{A}$.

$(2)$ is bicomplete if it is complete and cocomplete.

$(3)$ satisfies $\rm{AB4}$ if it is cocomplete and coproducts of short exact sequences (indexed over sets) in $\mathcal{A}$ are exact.
Dually, an abelian category $\mathcal{A}$ satisfies $\rm{AB4'}$ if it is complete
and products of short exact sequences (indexed over sets) in $\mathcal{A}$ are exact.
\end{Def}

If $\mathcal{A}$ satisfies $\rm{AB4}$, then $\D{\mathcal{A}}$ has coproducts indexed over sets, and therefore the coproducts of distinguished triangles in $\D{\mathcal{A}}$ are distinguished triangles. Moreover, $\D{\mathcal{A}}$ itself is the smallest full triangulated subcategory of $\D{\mathcal{A}}$ containing $\mathcal{A}$ and being closed under coproducts.
Dually, if $\mathcal{A}$ satisfies $\rm{AB4'}$, then $\D{\mathcal{A}}$ has products indexed over sets, and therefore the products of distinguished triangles in $\D{\mathcal{A}}$ are distinguished triangles.

Note that if a cocomplete (that is, coproducts indexed over sets exist) abelian category has enough injectives, then it satisfies $\rm{AB4}$. Dually, if a complete abelian category has enough projective, then it satisfies $\rm{AB4'}$. Examples of abelian categories with both $\rm{AB4}$ and $\rm{AB4'}$ are the module categories of rings, and the categories of additive functors from essentially small triangulated categories to the category of abelian groups  (see \cite[Chapter 6]{neemanbook} for details).

\begin{Lem}\label{HTP}
Let $\mathcal{A}$ be an abelian category and let $\mathcal{X}$ and $\mathcal{Y}$ be abelian subcategories of $\mathcal{A}$.
Assume that

$(1)$ $\mathcal{X}$ satisfies $\rm{AB4}$ and the inclusion $i: \mathcal{X}\subseteq \mathcal{A}$ preserves coproducts,

$(2)$ $\mathcal{Y}$ satisfies $\rm{AB4'}$ and the inclusion $j: \mathcal{Y}\subseteq \mathcal{A}$ preserves products, and

$(3)$ $\Hom_{\D{\mathcal{A}}}(X, Y[n])=0$ for all $X\in \mathcal{X}$, $Y\in\mathcal{Y}$ and $n\in\mathbb{Z}.$

\smallskip
\noindent Then $\Hom_{\D{\mathcal{A}}}(i(\cpx{X}), j(\cpx{Y}))=0$ for all $\cpx{X}\in\D{\mathcal{X}}$ and $\cpx{Y}\in\D{\mathcal{Y}}$.
\end{Lem}

{\it Proof.}  For any $Y\in\mathcal{Y}$, let $\mathscr{X}(Y)$ be the full subcategory of $\D{\mathcal{X}}$ consisting of objects $\cpx{X}$ such that
$\Hom_{\D{\mathcal{A}}}(i(\cpx{X}), j(Y)[n])=0$ for all $n\in\mathbb{Z}$. Then $\mathscr{X}(Y)$ is a full triangulated subcategory of $\D{\mathcal{X}}$.
By $(1)$, $\D{\mathcal{X}}$ has coproducts and $\mathscr{X}(Y)\subseteq\D{\mathcal X}$ is closed under coproducts. Moreover, $\mathcal{X}\subseteq \mathscr{X}(Y)$ by $(3)$. Note that $\D{\mathcal{X}}$ is the smallest full triangulated subcategory of $\D{\mathcal{X}}$ containing $\mathcal{X}$ and being closed under coproducts. Thus $\mathscr{X}(Y)=\D{\mathcal{X}}$. It follows that $\Hom_{\D{\mathcal{A}}}(i(\cpx{X}), j(Y)[n])=0$ for all $\cpx{X}\in\D{\mathcal{X}}$ and $n\in\mathbb{Z}$.
Dually, when $\cpx{X}\in\D{\mathcal{X}}$ is fixed, one can apply $(2)$ and $(3)$ to prove $\Hom_{\D{\mathcal{A}}}(i(\cpx{X}), j(\cpx{Y}))=0$ for all $\cpx{Y}\in\D{\mathcal{Y}}$. $\square$

\begin{Prop}\label{AB-condition}
 Suppose that $\mathcal{A}$ is an abelian category satisfying $\rm{AB4}$ and $\rm{AB4'}$. Let $\mathcal{X}$ and $\mathcal{Y}$ be abelian subcategories of $\mathcal{A}$ and let $*\in \{b, +,-, \varnothing\}$. Then the following are equivalent:

$(1)$ $(\mathcal{X}, \mathcal{Y})$ is a $\mathscr{D}^*$-decomposition of $\mathcal{A}$.

$(2)$ The inclusions $i:\mathcal{X}\subseteq\mathcal{A}$ and $j:\mathcal{Y}\subseteq\mathcal{A}$ induce fully faithful functors $\Ds{i}:\Ds{\mathcal{X}}\to \Ds{\mathcal{A}}$ and $\Ds{j}: \Ds{\mathcal{Y}}\to\Ds{\mathcal{A}}$, respectively, and $\big(\Img(\Ds{i}), \Img(\Ds{j})\big)$ is a semi-orthogonal decomposition of $\mathscr{D}^*(\mathcal{A})$.

\end{Prop}

{\it Proof.} $(2)$ implies $(1)$ by Lemma \ref{DD-HTP}.  Suppose $(1)$ holds. Then $(\mathcal{X}, \mathcal{Y})$ is a complete Ext-orthogonal pair in $\mathcal{A}$ by Proposition \ref{derived}. In particular, $\mathcal{X}={^\bot}\mathcal{Y}$. It follows that $\mathcal{X}$ is closed under coproducts in $\mathcal{A}$ and $i:\mathcal{X}\subseteq\mathcal{A}$ preserves coproducts. Since $\mathcal{A}$ satisfies $\rm{AB4}$, $\mathcal{X}$ also satisfies $\rm{AB4}$. Dually, $\mathcal{Y}$ is closed under products in $\mathcal{A}$ and satisfies $\rm{AB4'}$. Now, $(2)$ holds by Definitions \ref{Intro-def} and \ref{SOD} and Proposition \ref{HTP}. $\square$

\medskip
As a consequence of Theorem \ref{main-result}(2) and Proposition \ref{AB-condition}, we construct
half recollements of different derived categories from derived decompositions.

\begin{Koro}\label{half-recollement}
Let $\mathcal{A}$ be a bicomplete abelian category with enough projectives and injectives.
If $(\mathcal{X},\mathcal{Y})$ is a derived decomposition of $\mathcal{A}$, then there is a lower half recollement
$$\xymatrix{\Ds{\mathcal{X}}\ar^-{\Ds{i}}@/^0.8pc/[r]&\Ds{\mathcal{A}}
\ar^-{\mathbb{R}^*(r)}@/^0.8pc/[l]\ar^-{\mathbb{L}^*(\ell)}@/^0.8pc/[r]
&\Ds{\mathcal{Y}}\ar^-{\Ds{j}}@/^0.8pc/[l]}$$
for any $*\in \{b, +,-, \varnothing\}$, where $\mathbb{R}^*(r)$ and $\mathbb{L}^*(\ell)$ denote the right- and left-derived functors of $r$ and $\ell$, respectively.
\end{Koro}

\begin{Rem}
Consider the following statements:

$(1)$ $(\mathcal{X}, \mathcal{Y})$ is a $\mathscr{D}$-decomposition of $\mathcal{A}$;

$(2)$ $(\mathcal{X}, \mathcal{Y})$ is a $\mathscr{D}^*$-decomposition of $\mathcal{A}$ for any $\ast\in \{+,-\}$;

$(3)$ $(\mathcal{X}, \mathcal{Y})$ is a derived decomposition of $\mathcal{A}$.

\noindent Then $(1)\Rightarrow (2)\Rightarrow (3)$. Thus the existence of $\mathscr{D}^b$-decompositions is the weakest condition among those other type of derived decompositions introduced in Definition \ref{Intro-def}. This is why we sometimes pay more attention to the existence of such decompositions

To show the above implications, we consider the triangle given in Definition \ref{Intro-def}(D3). If $H^n(\cpx{M})=H^{n+1}(\cpx{M})=0$ for some integer $n$, then $H^{n+1}(X_{\cpx{M}})\simeq H^n(Y^{\cpx{M}})$. Clearly, $H^{n+1}(X_{\cpx{M}})\in\mathcal{X}$ and $H^n(Y^{\cpx{M}})\in\mathcal{Y}$ since $\mathcal{X}$ and $\mathcal{Y}$ are abelian subcategories of $\mathcal{A}$. However, $\mathcal{X}\cap\mathcal{Y}=\{0\}$ by Definition \ref{Intro-def}(D2). Thus  $H^{n+1}(X_{\cpx{M}})=H^n(Y^{\cpx{M}})=0$.

It is open whether $(3)$ always implies $(1)$. But Theorem \ref{main-result}(2) tells us this is true if $\mathcal{A}$ has enough projectives and injectives.
\end{Rem}

\subsection{One-to-one correspondence between derived and homological decompositions\label{3.4} }
In this section we establish a one-to-one correspondence between derived decompositions of abelian categories and homological decompositions, a subclass of semi-orthogonal decompositions, of bounded derived categories of abelian categories.

Let $\mathcal{A}$ be an abelian category. Given a full subcategory $\mathscr{X}$ of $\Db{\mathcal{A}}$, we denote by $\mathcal{X}_0$
the full subcategory of $\mathcal{A}$ consisting of $0$-th cohomologies of all objects in $\mathscr{X}$. Conversely, given a full subcategory $\mathcal{X}$
of $\mathcal{A}$, we denote by $\mathscr{D}^b_\mathcal{X}(\mathcal{A})$  the full subcategory of $\Db{\mathcal{A}}$ consisting of complexes such that all of their
cohomologies are in $\mathcal{X}$.

\begin{Def}\label{HD}
Let $\mathcal{A}$ be an abelian category. A pair $(\mathscr{X}, \mathscr{Y})$ of full subcategories of $\Db{\mathcal{A}}$ is called a \emph{homological decomposition} if the following three conditions hold:

$(a)$ $(\mathscr{X}, \mathscr{Y})$ is a semi-orthogonal decomposition of $\Db{\mathcal{A}}$.

$(b)$ $\mathcal{X}_0\subseteq \mathscr{X}$ and the inclusion $\mathcal{X}_0\subseteq\mathcal{A}$ induces a fully faithful functor $\Db{\mathcal{X}_0}\to \Db{\mathcal{A}}$.

$(c)$ $\mathcal{Y}_0\subseteq\mathscr{Y}$ and the inclusion $\mathcal{Y}_0\subseteq\mathcal{A}$ induces a fully faithful functor $\Db{\mathcal{Y}_0}\to \Db{\mathcal{A}}$.
\end{Def}
\begin{Lem}\label{GHOM}
Let $\mathscr{X}$ be a full triangulated subcategory of $\Db{\mathcal{A}}$. If $\mathcal{X}_0\subseteq\mathscr{X}$, then $\mathcal{X}_0$ is an abelian subcategory of $\mathcal{A}$ and $\mathscr{X}=\mathscr{D}^b_{\mathcal{X}_0}(\mathcal{A})$.
\end{Lem}

{\it Proof.} Since $\mathscr{X}$ is a full triangulated subcategory of $\Db{\mathcal{A}}$, it follows that  $H^n(\cpx{X})=H^0(\cpx{X}[n])\in\mathcal{X}_0$ for any $n\in\mathbb{Z}$ and $\cpx{X}\in\mathscr{X}$.  Suppose $\mathcal{X}_0\subseteq\mathscr{X}$. Then $H^n(\cpx{X})\in\mathscr{X}$. Let $f:X^{-1}\to X^0$ be a morphism in $\mathcal{A}$ with $X^{-1}, X^0\in\mathcal{X}_0$. Now, regarding $f$ as a complex in degrees $-1$ and $0$, we obtain $f\in \mathscr{X}$.
Consequently, as homologies of $f$, $\Ker(f)$ and $\Coker(f)$ belong to $\mathcal{X}_0$. Thus $\mathcal{X}_0\subseteq \mathcal{A}$ is an abelian subcategory.

Note that $\mathscr{X}\subseteq\mathscr{D}^b_{\mathcal{X}_0}(\mathcal{A})$, while the converse inclusion is due to the fact that each object $\cpx{M}$ in $\Db{\mathcal{A}}$ belongs to the smallest full triangulated subcategory containing $H^n(\cpx{M})$ for all $n$. $\square$

\begin{Prop}\label{Corre}
\noindent There is a one-to-one correspondence
$$\{\mbox{derived decompositions of}\; \mathcal{A}\}\longleftrightarrow \{\mbox{homological decompositions of}\;\Db{\mathcal{A}}\},$$
$$(\mathcal{X},\mathcal{Y})\mapsto(\mathscr{D}^b_\mathcal{X}(\mathcal{A}), \mathscr{D}^b_\mathcal{Y}(\mathcal{A})),$$
$$(\mathcal{X}_0, \mathcal{Y}_0)\leftmapsto (\mathscr{X}, \mathscr{Y}),$$
where $(\mathscr{X}, \mathscr{Y})$ is a homological decomposition of $\Db{\mathcal{A}}$ and $(\mathcal{X},\mathcal{Y})$ is a derived decomposition
of $\mathcal{A}$.
\end{Prop}

{\it Proof.} Given a homological decomposition $(\mathscr{X}, \mathscr{Y})$ of $\Db{\mathcal{A}}$, it follows from Definition \ref{HD} and Lemma \ref{GHOM} that $\mathscr{X}=\mathscr{D}^b_{\mathcal{X}_0}(\mathcal{A})$ and $\mathscr{Y}=\mathscr{D}^b_{\mathcal{Y}_0}(\mathcal{A})$.
Since the inclusion $u:\mathcal{X}_0\subseteq\mathcal{A}$ induces a fully faithful functor $\Db{\mathcal{X}_0}\to \Db{\mathcal{A}}$,
we have $\Img(\Db{u})=\mathscr{D}^b_{\mathcal{X}_0}(\mathcal{A})$ by Lemma \ref{subcat}. Similarly, $\Img(\Db{v})=\mathscr{D}^b_{\mathcal{Y}_0}(\mathcal{A})$
where $v:\mathcal{Y}_0\to\mathcal{A}$ is the inclusion. Thus $(\Img(\Db{u}), \Img(\Db{v}))$ is a semi-orthogonal decomposition of $\Db{\mathcal{A}}$.
By definition \ref{fdad}, $(\mathcal{X}_0, \mathcal{Y}_0)$ is a derived decomposition of $\mathcal{A}$.

Given a derived decomoposition $(\mathcal{X},\mathcal{Y})$ of $\mathcal{A}$, by Definition \ref{fdad} and Lemma \ref{subcat}, the pair $(\mathscr{D}^b_\mathcal{X}(\mathcal{A}), \mathscr{D}^b_\mathcal{Y}(\mathcal{A}))$ is a semi-orthogonal decomposition of $\Db{\mathcal{A}}$. Clearly, $(\mathscr{D}^b_\mathcal{X}(\mathcal{A}))_0=\mathcal{X}\subseteq \mathscr{D}^b_\mathcal{X}(\mathcal{A})$ and $(\mathscr{D}^b_\mathcal{Y}(\mathcal{A}))_0=\mathcal{Y}\subseteq \mathscr{D}^b_\mathcal{Y}(\mathcal{A})$.
By Definitions \ref{HD} and \ref{fdad}, $(\mathscr{D}^b_\mathcal{X}(\mathcal{A}), \mathscr{D}^b_\mathcal{Y}(\mathcal{A}))$ is a homological decomposition of $\Db{\mathcal{A}}$. Thus the correspondence in Proposition \ref{Corre} is well defined.
It is easy to see from Lemma \ref{GHOM} that the correspondence is one-to-one. $\square$

\section{Constructing derived decompositions\label{sect4}}
In this section we apply Theorem \ref{main-result} to construct derived decompositions for the module categories of rings. We first show that homological ring
epimorphisms can provide derived decompositions (see Proposition \ref{pd1}), and then prove that localizing subcategories and right perpendicular
subcategories in abelian categories also give rise to derived decompositions (see Proposition \ref{Localizing}). Finally, we construct derived decompositions
for module categories over left nonsingular rings and commutative noetherian rings (see Corollaries \ref{DADS} and \ref{Localizing commutative}, respectively). Moreover, this  construction establishes a derived stratification of
the module category of a commutative ring with the Krull dimension at most one (see Corollary \ref{decomposition}).

\subsection{Homological ring epimorphisms}
In this section we show that homological ring epimorphisms produce not only derived decompositions, but  also derived equivalences and recollements (see Corollary \ref{derived eqivalence}).

Throughout this section, we assume that $\lambda:R\to S$ is a homological ring epimorphism.
Define $$\begin{array}{rl}
\mathcal{A}:= & R\Modcat,\qquad \mathcal{X}:=  S\Modcat, \\
\mathcal{Y}:= & \{Y\in R\Modcat\mid\Hom_R(S, Y)=0=\Ext_R^1(S, Y)\}, \\ \mathcal{Z}:= &\{Z\in R\Modcat\mid S\otimes_RZ=0=\Tor_1^R(S, Z)\}.\end{array}$$

\begin{Prop}\label{pd1}
$(1)$ $(\mathcal{X}, \mathcal{Y})$ is a complete Ext-orthogonal pair in $\mathcal{A}$ if and only if $\pd(_RS)\leq 1$.

$(2)$  $(\mathcal{X}, \mathcal{Y})$ is a derived decomposition of $\mathcal{A}$
if and only if $\pd(_RS)\leq 1$ and $\Hom_R(\Coker(\lambda)$, $\Ker(\lambda))$ = $0.$
\end{Prop}

{\it Proof.} $(1)$ Since $\lambda$ is a ring epimorphism, the restriction functor $\lambda_*:\mathcal{X}\to\mathcal{A}$ is fully faithful. So, we identify $\mathcal{X}$ with the image of $\lambda_*$. Further, since $\lambda$ is homological, the derived functor $\Db{\lambda_*}:\Db{S}\to \Db{R}$ is fully faithful.
Note that ${_S}S\in \mathscr{P}(\mathcal{X})$, the category of projective $S$-modules. If $(\mathcal{X}, \mathcal{Y})$ is a complete Ext-orthogonal pair in $\mathcal{A}$, then $\pd({_R}S)\leq 1$
by Corollary \ref{proj or injective}(2). This shows the necessity of $(1)$.

To show the sufficiency of $(1)$, we assume $\pd(_RS)\leq 1$. Then $\mathcal{Y}=S^{\bot}$. It follows from \cite[Proposition 1.1]{GL} that $\mathcal{Y}$ is an abelian full subcategory of $\mathcal{A}$. Since ${^\bot}\mathcal{Y}$ contains ${_R}S$ and is closed under direct sums in $\mathcal{A}$, it must contain all projective $S$-modules. Moreover, each object of $\mathcal{X}$ admits a projective resolution by projective $S$-modules. Consequently, for any $X\in\mathcal{X}$, $Y\in\mathcal{Y}$ and $n\in\mathbb{N}$,  $\Hom_{\Db{R}}(X, Y[n])\simeq
\Hom_{\Db{R}}(\Omega_S^n(X), Y)=0$, where $\Omega_S^n(X)$ denotes an $n$-th syzygy module of ${_S}X$. This implies $\mathcal{X}\subseteq {^\bot}\mathcal{Y}$. By Lemma \ref{abelian}(2), to show $(1)$, it suffices to prove that $(\mathcal{X}, \mathcal{Y})$ satisfies $(GC)$.

The functor $\Db{\lambda_*}:\Db{S}\to \Db{R}$ has a right adjoint functor $\rHom_R(S, -):\Db{R}\to\Db{S}$. Let $\varepsilon: \Db{\lambda_*}\,\rHom_R(S, -)\ra \mbox{Id}_{\Db{R}}$ be the counit adjunction. Then, for each $R$-module $M$, there exists a distinguished triangle in $\D{R}$:
$$
(\dag)\quad \Db{\lambda_*}\,\rHom_R(S, M)\lraf{\varepsilon_M} M\lra \cpx{Y}_M\lra \Db{\lambda_*}\,\rHom_R(S, M)[1].
$$
Since $\lambda$ is homological, $\Db{\lambda_*}$ is fully faithful. So the morphism $\rHom_R(S,\varepsilon_M )$ is an isomorphism  in $\Db{S}$. This means $\cpx{Y}_M\in \mathscr{Y}:=\Ker(\rHom_R(S, -))\subseteq \Db{R}$. Thus
$$
\cpx{Y}_M\in\mathscr{Y}=\{\cpx{Y}\in\Db{R}\mid \Hom_{\Db{R}}(S, \cpx{Y}[n])=0 \;\;\mbox{for all}\; n\in\mathbb{Z}\}\quad \mbox{and}\quad \mathcal{Y}=\mathscr{Y}\cap R\Modcat.
$$
Taking cohomologies on the triangle ($\dag$) yields an exact sequence of $R$-modules:
$$
0\lra H^{-1}(\cpx{Y}_M)\lra \Hom_R(S,M)\lraf{\Hom_R(\lambda, M)} M\lra H^0(\cpx{Y}_M)\lra \Ext_R^1(S, M)\lra 0
$$
where $M$ is identified with $\Hom_R(R, M)$. Clearly, both $\Hom_R(S,M)$ and $\Ext_R^1(S,M)$ belong to $\mathcal{X}$. On the other hand, since $\pd(_RS)\leq 1$, the $R$-module $S$ is isomorphic in $\Db{R}$ to a two-term complex of projective $R$-modules and there is the following exact sequence by \cite[Lemma 3.4]{xc1}:
$$0\lra\Hom_{\Db{R}}(S,H^{n-1}(\cpx{Y})[1])\lra\Hom_{\Db{R}}(S,\cpx{Y}[n])\lra \Hom_{\Db{R}}(S,H^n(\cpx{Y}))\lra 0.$$
This shows $\mathscr{Y}=\{\cpx{Y}\in\Db{R}\mid H^n(\cpx{Y})\in\mathcal{Y}\;\;\mbox{for all}\; n\in\mathbb{Z}\}$. It then follows from  $\cpx{Y}_M\in\mathscr{Y}$ that $H^i(\cpx{Y}_M)\in\mathcal{Y}$ for any $i\in\mathbb{Z}$. Now, we define $X_M:=\Hom_R(S,M)$, $X^M:=\Ext_R^1(S, M)$, $Y_M:=H^{-1}(\cpx{Y}_M)$ and $Y^M:=H^0(\cpx{Y}_M)$. This shows the sufficiency of $(1)$.

$(2)$ Clearly, $\mathcal{A}$ has enough projectives and injectives. If $M$ is injective, then $X^M=0$. Note that $Y_M\simeq \Ker(\Hom_R(\lambda, M))\simeq \Hom_R(\Coker(\lambda), M)$ as $R$-modules. By $(1)$ and Theorem \ref{main-result}, $(\mathcal{X}, \mathcal{Y})$ is a derived decomposition of $\mathcal{A}$ if and only if $\pd({_R S})\leq 1$ and $\Hom_R(\Coker(\lambda), M)=0$ whenever $M$ is projective. To check $\Hom_R(\Coker(\lambda), M)=0$ for projective modules $M$, we only need to show $\Hom_R(\Coker(\lambda), R)=0$. because $\Hom_R(\Coker(\lambda), -)$ commutes with products and each projective $R$-module can be embedded into a product of copies of $R$. However, since $\lambda$ is a ring epimorphism,  $\Hom_R(\Coker(\lambda), S)=0$. This implies $\Hom_R(\Coker(\lambda),R)\simeq\Hom_R(\Coker(\lambda), \Ker(\lambda))$. Thus $(2)$ holds. $\square$

\medskip
When dealing with flat dimensions instead of projective dimensions, we have the following

\begin{Prop}\label{fld1}
$(1)$ $(\mathcal{Z}, \mathcal{X})$ is a complete Ext-orthogonal pair in $\mathcal{A}$ if and only if $\fld(S_R)\leq 1$.

$(2)$ $(\mathcal{Z}, \mathcal{X})$ is a derived decomposition of $\mathcal{A}$
if and only if $\fld(S_R)\leq 1$ and $\Coker(\lambda)\otimes_RI=0$ for any injective $R$-module $I$.
\end{Prop}

{\it Proof.} The proof of this result is similar to the one of Proposition \ref{pd1}. For the convenience of the reader, we list some key points in the proof.

Let $J:=\Hom_\mathbb{Z}(S_S, \mathbb{Q}/\mathbb{Z})$. Then $J$ is an injective cogenerator in $S\Modcat.$
Further, $\id({_R}J)=\fld(S_R)$ because a right $R$-module $N$ is flat
if and only if ${_R}\Hom_\mathbb{Z}(N, \mathbb{Q}/\mathbb{Z})$ is injective.
So the necessity of $(1)$ follows from Corollary \ref{proj or injective}(1).
If $\fld(S_R)\leq 1$, then the following statements hold true:

$(1)$  $\mathcal{Z}$ is abelian full subcategory of $\mathcal{A}$.

$(2)$  For any $\cpx{M}\in\D{R}$, $S\otimesL_R\cpx{M}=0$ if and only if $H^n(\cpx{M})\in\mathcal{Z}$ for all $n\in\mathbb{Z}$.

$(3)$ For any $M\in\mathcal{A}$, there is an exact sequence of $R$-modules:
$$
0\lra \Tor_1^R(S,M)\lra X_M\lra M\lraf{\lambda\otimes_RM} S\otimes_RM\lra X^M\lra 0
$$
such that $X_M, X^M\in\mathcal{Z}$. Clearly, $X^M\simeq \Coker(\lambda)\otimes_RM$ as $R$-modules.
Now, all other assertions in Proposition \ref{fld1} can be concluded from Theorem \ref{main-result}. $\square$

\medskip
As a consequence of Proposition \ref{fld1}, we have the following result on localizations of commutative rings.

\begin{Koro}\label{commtative}
Let $R$ be a commutative noetherian ring with $\Phi$ a multiplicative
subset of $R$, $\Phi^{-1}R$ the localization of $R$ at
$\Phi$ and $\mathcal{U}:=\{X\in R\Modcat\mid (\Phi^{-1}R)\otimes_RX=0\}$. Then $\big(\mathcal{U}, (\Phi^{-1}R)\Modcat\big)$ is a derived decomposition of $R\Modcat$.
\end{Koro}

{\it Proof.} Let $\lambda: R\to S:=\Phi^{-1}R$ be the localization of $R$ at
$\Phi$. Then $S$ is commutative and flat as an $R$-module, and therefore $\lambda$ is a homological ring epimorphism. By Proposition \ref{fld1}(2), it suffices to show  $\Coker(\lambda)\otimes_RI=0$ (or equivalently, $\lambda\otimes_RI$ is surjective) for any injective $R$-module $I$.

Since $R$ is a commutative noetherian ring, each injective $R$-module is a direct sum of indecomposable injective $R$-modules (see \cite[Theorem 3.3.10]{EJ}). So we only need to check the surjection of $\lambda\otimes_RI$  whenever $I$ is indecomposable.
By \cite[Theorem 3.3.7]{EJ}, there is a prime ideal $\mathfrak{p}$ of $R$ such that $I$ is isomorphic to the injective envelope $E(R/\mathfrak{p})$ of the $R$-module $R/\mathfrak{p}$. Moreover, by \cite[Theorem 3.3.8(6)]{EJ}, $\lambda\otimes_R E(R/\mathfrak{p})$ is an isomorphism if
$\Phi\cap \mathfrak{p}=\emptyset$; and $S\otimes E(R/\mathfrak{p})=0$ if $\Phi\cap \mathfrak{p}\ne\emptyset$. This implies that $\lambda\otimes_RE(R/\mathfrak{p})$ is always surjective, and therefore $\lambda\otimes_RI$ is surjective. Thus Corollary \ref{commtative} follows from Proposition \ref{fld1}(2). $\square$

\medskip
Next, we show that the derived decompositions in Propositions \ref{pd1} and \ref{fld1}
provide also lower half recollements of derived categories.

For the ring epimorphism $\lambda:R\to S$, we consider it as a complex $\cpx{Q}: 0\ra R\lraf{\lambda} S\ra 0$
of $R$-$R$-bimodules with $R$ and $S$ in degrees $-1$ and $0$, respectively.
Let
$$F:=\cpx{Q}[-1]\otimesL_R-:\D{R}\to\D{R},\quad G:=\rHom_R(\cpx{Q}[-1], -):\D{R}\to\D{R}$$
and let ${\rm Tria}(_R\cpx{Q})$ be the smallest full triangulated subcategory of $\D{R}$ containing $\cpx{Q}$ and being closed under direct sums.
Then $(F,G)$ is an adjoint pair of triangle functors and the restriction of $G$ to ${\rm Tria}(_R\cpx{Q})$ is fully faithful (see \cite[Section 4]{NS}).

In the case of Proposition \ref{pd1}(2), it follows from Corollary \ref{half-recollement} that there is a lower half recollement of bounded derived categories:
$$(\ddag)\quad
\xymatrix{\Ds{S}\ar^-{\Ds{\lambda_*}}@/^0.8pc/[r]&\Ds{R}
\ar^-{{\rm\mathbb{R}}^*\Hom_R(S,-)}@/^0.8pc/[l]\ar^-{\mathbb{L}^*(\ell)}@/^0.8pc/[r]
&\Ds{\mathcal{Y}}\ar^-{\Ds{j}}@/^0.8pc/[l]}
$$
for any $*\in \{b, +,-, \varnothing\}$, where $\ell:R\Modcat\to\mathcal{Y}$ is a left adjoint of the inclusion $j:\mathcal{Y}\to R\Modcat$. We claim that $\ell$ is the composition of the functors:
$$
R\Modcat\hookrightarrow\Ds{R}\lraf{G}\Ds{R}\lraf{H^0} R\Modcat.
$$
In fact, there is a canonical triangle
$\cpx{Q}[-1]\lraf{\sigma} R\lraf{\lambda} S\lraf{\pi} \cpx{Q}$ in
$\Db{R\otimes_\mathbb{Z}R\opp}$ which induces a sequence of triangle functors from $\Ds{R}$ to $\Ds{R}$
$$(\ddag)\quad
G[-1]\lraf{\pi_*} \D{\lambda_*}\rHom_R(S, -) \lraf{\lambda_*} \rHom_R(R, -)\lraf{\sigma_*} G
$$
such that their operations on a fixed object in $\Ds{R}$ yield a triangle in $\Ds{R}$. Clearly, $\rHom_R(R, -)$ can be identified with
the identity functor of $\Ds{R}$ up to natural isomorphism. So, for an $R$-module $M$, by taking cohomologies on ($\ddag$), we get a long exact sequence of $R$-modules:
$$
0\lra H^{-1}(G(M))\lra \Hom_R(S,M)\lraf{\Hom_R(\lambda, M)} M\lra H^0(G(M))\lra \Ext_R^1(S, M)\lra 0.
$$
Now, as in the proof of Proposition \ref{pd1}(1), both $\Hom_R(S,M)$ and $\Ext_R^1(S, M)$ belong to $\mathcal{X}$ and both $H^{-1}(G(M))$ and $H^0(G(M))$ belong to $\mathcal{Y}$. Thus $\ell(M)=H^0(G(M))$ by the definition of $\ell$.

Similarly, in the case of Proposition \ref{fld1}(2), we obtain a lower half recollement:
$$
\xymatrix{\Ds{\mathcal{Z}}\ar^-{\Ds{i}}@/^0.8pc/[r]&\Ds{R}
\ar^-{\mathbb{R}^*(r)}@/^0.8pc/[l]\ar^-{S\otimesL_R-}@/^0.8pc/[r]
&\Ds{S}\ar^-{\Ds{\lambda_*}}@/^0.8pc/[l]}
$$
where $r:=H^0 F(-):R\Modcat\to\mathcal{Z}$ is a right adjoint of the inclusion $i:\mathcal{Z}\to R\Modcat.$ The five-term exact sequence of $R$-modules is given by
$$
0\lra \Tor_1^R(S,M)\lra H^0(F(M))\lra M\lraf{\lambda\otimes_RM} S\otimes_RM\lra H^1(F(M)) \lra 0.
$$
for $M\in R\Modcat.$

The following result is an immediate consequence of Propositions \ref{pd1}(2) and \ref{fld1}(2), which supplies recollements and triangle equivalences of bounded derived categories.

\begin{Koro}\label{derived eqivalence}
Suppose

$(a)$ $\pd(_RS)\leq 1$ and $\Hom_R(\Coker(\lambda), \Ker(\lambda))=0$, and

$(b)$ $\fld(S_R)\leq 1$ and $\Coker(\lambda)\otimes_RI=0$ for any injective $R$-module $I$. Then

$(1)$ there is a recollement of $*$-bounded derived categories of abelian categories for $*\in \{b,+,-, \varnothing\}$:
$$\xymatrix@C=1.3cm{\Ds{S}\ar[r]^-{\Ds{\lambda_*}}&\Ds{R}\ar[r]^-{\mathbb{L}^*(\ell)}
\ar@/^1.2pc/[l]\ar@/_1.8pc/[l]
&\Ds{\mathcal{Y}}\ar@/^1.2pc/[l]\ar@/_1.8pc/[l]^-{}\vspace{0.3cm}}.$$

\medskip
$(2)$ The functor $\mathbb{R}^*(r)\Ds{j}: \Ds{\mathcal{Y}}\lraf{\simeq} \Ds{\mathcal{Z}}$
is a triangle equivalence with a quasi-inverse functor $\mathbb{L}^*(\ell)\Ds{i}$ for $*\in \{b,+,-, \varnothing\}$.
\end{Koro}

An example of Corollary \ref{derived eqivalence} is the following: Let $R$ be a $1$-Gorenstein ring (that is, a commutative noetherian ring such that the injective dimension of $R$ is at most $1$) and let $\Phi$ be the set of all non-zero divisors of $R$. Then $\lambda$ is always injective and ${_R}S$ is flat, injective and of projective dimension at most $1$. In case of the inclusion $\mathbb{Z}\subseteq \mathbb{Q}$, we get a recollement $(\Ds{\mathbb{Q}},\Ds{\mathbb{Z}},\Ds{\mathcal{Y}})$ and an equivalence $\Ds{\mathcal{Y}}\simeq \Ds{\mathcal{Z}}$.

\begin{Rem}
Corollary \ref{derived eqivalence}(2) deals with arbitrary homological ring epimorphisms. It generalizes the Matlis equivalences of
derived categories of abelian categories in \cite[Theorem 7.6]{Pos} for the localization $\lambda: R\ra S$ of a commutative ring $R$ at a multiplicative subset $\Phi$. In fact, if $\pd({_R}S)\leq 1$ and the $\Phi$-torsion in $R$ is bounded (that is, $\Ker(\lambda)$ is annihilated by an element in $\Phi$), then $(a)$ and $(b)$ in Corollary \ref{derived eqivalence} are satisfied by \cite[Lemma 6.1]{Pos}, and therefore $\Ds{\mathcal{Y}}\lraf{\simeq} \Ds{\mathcal{Z}}$.

Very recently, the derived equivalences
of bounded derived categories in Corollary \ref{derived eqivalence}(2) has been extended to the ones of derived categories of other types in \cite[Corollary 18.5]{BP}. Thus Corollary \ref{derived eqivalence}(2)
coincide with \cite[Corollary 18.5]{BP}.
\end{Rem}

\subsection{Localizing subcategories}\label{Section 4.2}
In this section we construct derived decompositions from localizing subcategories.

Let $\mathcal{A}$ be an abelian category and $\mathcal{X}$  a
full subcategory of $\mathcal{A}$. We say that $\mathcal{X}$ is a \emph {Serre} subcategory if it is closed under subobjects, quotients and extensions. In particular, $\mathcal{X}$ is an abelian
subcategory of $\mathcal{A}$, and the \emph {quotient category} $\mathcal{A}/\mathcal{X}$
(in the sense of Gabriel, Grothendieck, Serre) is defined by inverting all these morphisms in $\mathcal{A}$ that have kernels and cokernels
in $\mathcal{X}$. The quotient category has the same objects as $\mathcal{A}$ and is again an abelian category.
Moreover, there is a canonical exact functor $q:\mathcal{A}\to \mathcal{A}/\mathcal{X}$ (called the quotient functor)
such that the kernel of $q$ is exactly $\mathcal{X}$.

A Serre subcategory $\mathcal{X}$ of $\mathcal{A}$ is called a \emph{localizing} subcategory of $\mathcal{A}$ if $q$ has a right adjoint $s:\mathcal{A}/\mathcal{X}\to \mathcal{A}$
(called the section functor). This is equivalent to saying that $q$ restricts to an equivalence of additive categories
from $\mathcal{X}^{\bot {0,1}}:=\mathcal{X}^{\bot 0}\cap\mathcal{X}^{\bot 1}$ to $\mathcal{A}/\mathcal{X}$ (see \cite[Chap. III.2]{Gabriel} and \cite[Proposition 2.2]{GL}). In this case, $\mathcal{X}={^{\bot {0,1}}}(\mathcal{X}^{\bot {0,1}})$.
Note that $\mathcal{X}^{\bot {0,1}}$ is closed under extensions and kernels in $\mathcal{A}$ (see, for example, \cite[Proposition 1.1]{GL}), but it may not be an abelian subcategory of $\mathcal{A}$ in general.

If $\mathcal{A}$ is a Grothendieck category (that is, an ableian category with a generator and coproducts such that direct limits of exact sequences are exact), then a Serre subcategory of $\mathcal{A}$ is localizing if and only if it is
closed under coproducts in $\mathcal{A}$ (see \cite[Proposition 2.5]{GL}).

\begin{Lem}\label{Serre Case}
Let $\mathcal{A}$ be an abelian category and $\mathcal{X}$ a localizing subcategory of $\mathcal{A}$ with $\mathcal{Y}:=\mathcal{X}^{\bot}$.
Then

$(1)$ $(\mathcal{X}, \mathcal{Y})$ is a complete Ext-orthogonal pair in $\mathcal{A}$
if and only if $\mathcal{Y}=\mathcal{X}^{\bot {0,1}}$ if and only if  the section functor $s:\mathcal{A}/\mathcal{X}\to \mathcal{A}$ is exact.

$(2)$ $(\mathcal{X}, \mathcal{Y})$ is a derived decomposition of  $\mathcal{A}$
if and only if both $\mathcal{Y}=\mathcal{X}^{\bot {0,1}}$ and the derived functor $\Db{i}:\Db{\mathcal{X}}\to \Db{\mathcal{A}}$,
induced from the inclusion $i:\mathcal{X}\to\mathcal{A}$, is fully faithful.

\end{Lem}
{\it Proof.} $(1)$ Since $\mathcal{X}$ is a localizing subcategory of $\mathcal{A}$, it follows from
\cite[Proposition 2.2]{GL} that, for each object $M\in\mathcal{A}$, there is an exact sequence
$0\to X_1\to M\to \overline{M}\to X_2\to 0$ in $\mathcal{A}$ with $X_1, X_2\in\mathcal{X}$ and $\overline{M}\in\mathcal{X}^{\bot {0,1}}$.
Clearly, $\mathcal{X}\subseteq {}^\bot\mathcal{Y}$. By Lemma \ref{abelian}, $(\mathcal{X}, \mathcal{Y})$ is a complete Ext-orthogonal pair in $\mathcal{A}$ if and only if $\mathcal{Y}=\mathcal{X}^{\bot {0,1}}$. It remains to show that $\mathcal{Y}=\mathcal{X}^{\bot {0,1}}$ if and only if $s$ is exact.

Let $\mathcal{B}:=\mathcal{A}/\mathcal{X}$ and $q_1:\mathcal{X}^{\bot {0,1}}\lraf{\simeq}\mathcal{B}$ the restriction of the canonical functor $q: \mathcal{A}\to \mathcal{A}/\mathcal{X}$ to $\mathcal{X}^{\bot {0,1}}$. It is known that $s$ is always fully faithful and isomorphic to the composition of the quasi-inverse of $q_1$ with the inclusion $\mathcal{X}^{\bot {0,1}}\subseteq\mathcal{A}$ (see, for example, \cite[Proposition 2.2]{GL}).  If $\mathcal{Y}=\mathcal{X}^{\bot {0,1}}$, then $\mathcal{Y}$ is an abelian subcategory of $\mathcal{A}$ since $\mathcal{X}^{\bot {0,1}}$ is closed under extensions and kernels in $\mathcal{A}$. In this case, $q_1$ is an equivalence of abelian categories, and thus $s$ is exact.

Conversely, suppose that $s$ is an exact functor. Since both $q$ and $s$ are exact, they induce derived functors $\Db{q}:\Db{\mathcal{A}}\to\Db{\mathcal{B}}$ and $\Db{s}:\Db{\mathcal{B}}\to \Db{\mathcal{A}}$ such that $(\Db{q}, \Db{s})$ is an adjoint pair and $\Db{s}$ is fully faithful.
Picking up an object $Y\in \mathcal{X}^{\bot {0,1}}$, we then have $Y\simeq s(Z)$ for some $Z\in \mathcal{B}$ since $\Img(s)=\mathcal{X}^{\bot {0,1}}$.
For any $X\in\mathcal{X}$ and $n\in\mathbb{N}$,
$$\Ext_\mathcal{A}^n(X, Y)\simeq \Ext_\mathcal{A}^n(X, s(Z))=\Hom_{\Db{\mathcal A}}(X, s(Z)[n])\simeq\Hom_{\Db{\mathcal B}}(q(X), Z[n])=\Hom_{\Db{\mathcal B}}(0, Z[n])=0.$$ This implies both $\mathcal{X}^{\bot {0,1}}\subseteq \mathcal{Y}$ and $\mathcal{X}^{\bot {0,1}}=\mathcal{Y}$.

$(2)$ The necessity of the conditions in $(2)$  is a consequence of Definition \ref{Intro-def}, Proposition \ref{derived} and $(1)$.
Now, we show the sufficiency of the conditions in $(2)$.

Suppose that the functor $\Db{i}$ is fully faithful and $\mathcal{Y}=\mathcal{X}^{\bot {0,1}}.$ Let $\mathscr{X}:=\Ker(\Db{q})$ and let
$\eta:{\rm Id}_{\Db{\mathcal{A}}}$ $\to \Db{s}\Db{q}$ and $\epsilon: \Db{q}\Db{s}\to \mbox{Id}_{\Db{\mathcal{B}}}$ be the
unit and counit of the adjoint pair $(\Db{q}, \Db{s}),$ respectively. The full faithfulness of $\Db{s}$ is equivalent to saying that the morphism $\epsilon_{\cpx{Z}}: \Db{q}\Db{s}(\cpx{Z})\ra \cpx{Z}$ is an isomorphism for each $\cpx{Z}\in\Db{\mathcal{B}}$ (see \cite[Theorem 1, p.90]{Maclane}).
Now, for each $\cpx{L}\in\Db{\mathcal{A}}$, we take a triangle $\cpx{N}\to \cpx{L}\lraf{\eta_{\cpx{L}}} \Db{s}\Db{q}(\cpx{L})\to \cpx{N}[1]$ in $\Db{\mathcal{A}}$. It then follows from $\Db{q}(\eta_{\cpx{L}})\epsilon_{\Db{q}(\cpx{L})}=id_{\Db{q}\Db{s}\Db{q}(\cpx{L})}$ (see \cite[(9), p.85]{Maclane}) that $\Db{q}(\eta_{\cpx{L}})$ is an isomorphism and so $\cpx{N}\in\mathscr{X}$. Further, for any $\cpx{X}\in\mathscr{X}$ and $\cpx{Y}\in \Img(\Db{s})$, the adjunction isomorphism of the pair $(\Db{q}, \Db{s})$ implies $\Hom_{\Db{\mathcal{A}}}(\cpx{X},\cpx{Y})=0$. Consequently, $(\mathscr{X}, \Img(\Db{s}))$ is a semi-orthogonal decomposition of
$\Db{\mathcal{A}}$. Since $q$ is exact, $\mathscr{X}$ coincides with the full triangulated subcategory of $\Db{\mathcal{A}}$
consisting of complexes $\cpx{X}\in \D{\mathcal{A}}$ such that $H^n(\cpx{X})\in \mathcal{X}$ for all $n$. Further, since $\Db{i}$ is fully faithful, $\mathscr{X} =\Img(\Db{i})$ by Lemma \ref{subcat}. Recall that $s$ is fully faithful and $\Img(s)=\mathcal{Y}$. Thus
$(\mathcal{X}, \mathcal{Y})$ is a derived decomposition of  $\mathcal{A}$. $\square$

 \medskip
\begin{Koro}\label{expansion}
Let $\mathcal{A}$ be an abelian category and $\mathcal{X}$ be a localizing subcategory of $\mathcal{A}$ with $\mathcal{Y}:=\mathcal{X}^{\bot}$.
Assume $\Ext_\mathcal{A}^2(X, M)=0$ for $X\in\mathcal{X}$ and $M\in\mathcal{A}$. Then $(\mathcal{X}, \mathcal{Y})$ is a derived decomposition of  $\mathcal{A}$
\end{Koro}

{\it Proof.} We first prove $\Ext_\mathcal{A}^n(X, M)=0$ for all $n\geq 2$, $X\in\mathcal{X}$ and $M\in\mathcal{A}$. This particularly implies the equality $\mathcal{Y}=\mathcal{X}^{\bot {0,1}}$. For $n=2$, this is true by assumption. For $n\ge 3$, we need the following general result:

If $\cpx{M}\in\Cb{\mathcal{A}}$ satisfies both $H^m(\cpx{M})=0$ for all $m\neq 0$ and $H^0(\cpx{M})\in\mathcal{X}$,
then there is a two-term complex $\cpx{E}: 0\to E^{-1}\to E^0\to 0$ in $\Cb{\mathcal{A}}$ and a quasi-isomorphism from $\cpx{E}$ to $\cpx{M}$. Moreover, if, in addition, $\cpx{M}\in \Cb{\mathcal{X}}$, then $\cpx{E}\in \Cb{\mathcal{X}}$.

Indeed, let $\cpx{M}:=(M^m, d^m)_{m\in\mathbb{Z}}$ with $d^m:M^m\to M^{m+1}$ and $N:=H^0(\cpx{M})$. Then the
inclusion from $\tau_{\leq 0}(\cpx{M})$ into $\cpx{M}$ is a quasi-isomorphism, where the truncated complex
$\tau_{\leq 0}(\cpx{M})$ at degree $0$ is of the form $\cdots\to M^{-2}\lraf{d^{-2}} M^{-1}\lraf{d^{-1}}\Ker(d^0)\to 0$. Since $\Ext_\mathcal{A}^2(N, \Img(d^{-2}))=0$,
we apply $\Ext_\mathcal{A}^1(N,-)$ to the exact sequence $0\to \Img(d^{-2})\to M^{-1}\to\Img(d^{-1})\to 0$ and obtain a surjective map
$\Ext_\mathcal{A}^1(N, d^{-1}): \Ext_\mathcal{A}^1(N, M^{-1})\to\Ext_\mathcal{A}^1(N, \Img(d^{-1}))$.
As the exact sequence $\delta:0\to \Img(d^{-1})\to \Ker(d^0)\to N\to 0$ corresponds to an element in
$\Ext_\mathcal{A}^1(N, \Img(d^{-1}))$, there is another exact sequence $\delta':0\to M^{-1}\lraf{g} E^0\to  N\to 0$ which is sent
to $\delta$ by taking the push-out of $(g, d^{-1})$. This also provides an associated morphism $h:E^0\to \Ker(d^0)$ such that
$gh=d^{-1}:M^{-1}\to\Ker(d^0)$. Let $\cpx{E}: 0\to E^{-1}\lraf{g} E^0\to 0$ with $E^{-1}:=M^{-1}$, and let $\cpx{f}:=(f^m)_{m\in\mathbb{Z}}$ with
$f^{-1}={\rm Id}_{M^{-1}}$, $f^0=h$ and $f^m=0$ for $m\neq 0, 1$. Then $\cpx{f}: \cpx{E}\to\tau_{\leq 0}(\cpx{M})$ is a quasi-isomorphism, and therefore the composition of $\cpx{f}$ with the inclusion $\tau_{\leq 0}(\cpx{M})\to\cpx{M}$ is a quasi-isomorphism. Since $\mathcal{X}$ is an abelian full subcategory
of $\mathcal{A}$ and is closed under extensions, we have $\cpx{E}\in \Cb{\mathcal{X}}$ if $\cpx{M}\in \Cb{\mathcal{X}}$.

Recall that $\Ext_\mathcal{A}^n(X, M)$ stands for $\Hom_{\Db{\mathcal{A}}}(X, M[n])$. By calculating this Hom-space
in $\Db{\mathcal{A}}$, we conclude from the above general result that $\Ext_\mathcal{A}^n(X, -)=0=\Ext_\mathcal{X}^n(X, -)$ for all $n\geq 2$
and $X\in\mathcal{X}$.

Thus the functor $\Db{i}:\Db{\mathcal{X}}\to\Db{\mathcal{A}}$, induced from the inclusion $i:\mathcal{X}\to\mathcal{A}$,
preserves $\Ext^n$ for all $n\geq 2$. It also preserves $\Hom$ and $\Ext^1$ because $\mathcal{X}$ is an abelian full subcategory
of $\mathcal{A}$ and is closed under extensions. Hence $\Db{i}$ induces $\Ext_\mathcal{X}^n(X, X')\lraf{\simeq} \Ext_\mathcal{A}^n(X, X')$
for all $n\geq 0$ and $X, X'\in\mathcal{X}$. This implies that $\Db{i}$ is fully faithful by Lemma \ref{embedding}.
Now, Corollary \ref{expansion} is a consequence of Lemma \ref{Serre Case}(2). $\square$

\begin{Lem}\label{Compare}
Suppose that $\mathcal{A}$ is an abelian category such that each object of $\mathcal{A}$ has an injective envelope.
Let $\mathcal{X}$ be a Serre subcategory of $\mathcal{A}$. Define  $\mathcal{E}:=\mathcal{X}^{\bot 0}\cap \mathscr{I}(\mathcal{A})$.
Then $\mathcal{X}^{\bot}$ (respectively, $\mathcal{X}^{\bot {0,1}}$) consists of all objects $M$ which has a minimal injective resolution
$0\to M\to I_0\to I_1\to \cdots \to I_i\to \cdots$
with $I_i\in\mathcal{E}$ for all $i\geq 0$ (respectively, $i=0,1$).
\end{Lem}

{\it Proof.} We first prove that $\mathcal{X}^{\bot 0}$ is closed under injective envelope in $\mathcal{A}$, that is,
if $Z\in \mathcal{X}^{\bot 0}$, then the injective envelope $E(Z)$ of $Z$ belongs to $\mathcal{X}^{\bot 0}$.

Let $Z\in\mathcal{X}^{\bot 0}$ and assume contrarily that there is a nonzero morphism $f: X\to E(Z)$ in $\mathcal{A}$ for some $X\in\mathcal{X}$.
Then $\Img(f)\ne 0$ and there is a monomorphism $g:\Img(f)\to E(Z)$. Let $h:Z\to E(Z)$ be an injective envelope of $Z$. Taking the pull-back of
$(g,h)$ yields another two monomorphisms $K\to Z$ and $K\to \Img(f)$ in $\mathcal{A}$. As $E(Z)$ is the injective envelope of $Z$, we have $K\ne 0$. By assumption, $\mathcal{X}$ is closed under subobjects and quotients. Hence, with $X$ also $\Img(f)$ and $K$ lie in $\mathcal{X}$. It follows
from $Z\in\mathcal{X}^{\bot 0}$ that $K=0$, a contradiction. This shows $E(Z)\in\mathcal{X}^{\bot 0}$. Hence $\mathcal{X}^{\bot 0}$ is closed under injective envelope in $\mathcal{A}$.

If $Z\in \mathcal{X}^{\bot 0}$, then $E(Z)\in\mathcal{E}$. Moreover, there are inclusions of categories: $\mathcal{E}\subseteq\mathcal{X}^{\bot}\subseteq \mathcal{X}^{\bot {0,1}}\subseteq\mathcal{A}$.
Recall that $\mathcal{X}^{\bot}$ is closed under extensions, kernels of epimorphisms and cokernels of monomorphisms
in $\mathcal{A}$, and that $\mathcal{X}^{\bot {0,1}}$ is closed under extensions and kernels in $\mathcal{A}$.
Now, it is easy to verify Lemma \ref{Compare}. $\square$

\medskip
The following result furnishes a way to get derived decompositions from localizing subcategories.

\begin{Prop}\label{Localizing}
Let $\mathcal{A}$ be an abelian category such that each of its objects has an injective envelope. If $\mathcal{X}$ is a localizing
subcategory of $\mathcal{A}$ with $\mathcal{Y}:=\mathcal{X}^{\bot}$, then the following are equivalent:

$(1)$ $(\mathcal{X}, \mathcal{Y})$ is a derived decomposition of $\mathcal{A}$.

$(2)$ Each morphism $I^0\to I^1$ between injective objects in $\mathcal{A}$ with $I^1\in\mathcal{Y}$
can be completed to an exact sequence $I^0\to I^1\to I^2$ such that $I^2$ is injective and $I^2\in\mathcal{Y}$.

$(3)$ The image of each morphism from an injective object in $\mathcal{A}$ to an object in $\mathcal{Y}$
belongs to $\mathcal{Y}$.
\end{Prop}

{\it Proof.}
Since $\mathcal{X}$ is a localizing subcategory of $\mathcal{A}$, the proof of Lemma \ref{Serre Case}(1)
shows that the five-term exact sequence associated with an object $M\in\mathcal{A}$ becomes
$$0\lra r(M)\lraf{\varepsilon_M^{-1}} M
\lraf{\varepsilon_M^{0}} \ell(M)\lra X^M\lra 0,$$
with $r(M),X^M\in \mathcal{X}$ and $\ell(M)\in \mathcal{X}^{\bot 0,1}$.

$(1)\Rightarrow (2)$: Since $\mathcal{A}$ has enough injectives and $\Db{i}:\Db{\mathcal{X}}\to \Db{\mathcal{A}}$ is fully faithful,
we see from Lemma \ref{EOff}(2) that $X^M=0$ whenever $M\in\mathscr{I}(\mathcal{A})$. Let $f:I^0\to I^1$ be a morphism of injective objects
in $\mathcal{A}$ with $I^1\in\mathcal{Y}$. Then there is an exact sequence $0\to r(I^0)\to I^0\to \ell(I^0)\to 0$ in $\mathcal{A}$.
We always have $\Hom_\mathcal{A}(r(I^0), I^1)=0$, due to $r(I^0)\in\mathcal{X}$ and $I^1\in\mathcal{Y}$.
Consequently, $f$ is the composition of the morphism $\varepsilon_{I^0}^{0}:I^0\to \ell(I^0)$
with another morphism $g: \ell(I^0)\to I^1$. This implies $\Coker(f)\simeq \Coker(g)$. Since $\mathcal{Y}$ is an abelian subcategory of $\mathcal{A}$,
$\Coker(g)\in\mathcal{Y}$, and thus also $\Coker(f)\in\mathcal{Y}$. Let $I^2$ be the injective envelope of $\Coker(f)$.
Then $I^2\in\mathcal{Y}$ by Lemma \ref{Compare}. Now, we extend $f$ to an exact sequence $I^0\to I^1\to I^2$.

$(2)\Rightarrow (3)$:
Thanks to Lemma \ref{Compare}, the assumption $(2)$ implies $\mathcal{Y}=\mathcal{X}^{\bot {0,1}}$.
Since $\mathcal{Y}$ is closed under extensions, kernels of epimorphisms and cokernels of monomorphisms in $\mathcal{A}$ and since $\mathcal{X}^{\bot {0,1}}$ is closed under kernels in $\mathcal{A}$ by \cite[Proposition 1.1]{GL}, $\mathcal{Y}$ is an abelian subcategory of $\mathcal{A}$.
Let $h:I^0\to Y$ be morphism in $\mathcal{A}$ with $Y\in\mathcal{Y}$ and $I^0$ an injective object.
Further, let $I^1$ be the injective envelope of $Y$ with a monomorphism $s:Y\to I^1$. Then $I^1\in\mathcal{Y}$, according to Lemma \ref{Compare}.
Moreover, by $(2)$, the composition $h$ with $s$ can be completed to an exact sequence $I^0\lraf{hs} I^1\lraf{t} I^2$ in $\mathcal{A}$ such that $I^2$ is injective and
$I^2\in\mathcal{Y}$. Thus $\Img(hs)=\Ker(t)\in\mathcal{Y}$. Since $\Img(h)\simeq\Img(hs)$,  $\Img(h)\in\mathcal{Y}$, and therefore $(3)$ follows.

$(3)\Rightarrow (1)$:
By Lemma \ref{Compare}, $\mathcal{Y}=\mathcal{X}^{\bot {0,1}}$. Thus, by Lemma \ref{Serre Case}(2), to show $(1)$, it suffices to prove that $\Db{i}:\Db{\mathcal{X}}\to \Db{\mathcal{A}}$ is fully faithful.

In fact, given an injective object $I$ of $\mathcal{A}$, since $\ell(I)\in\mathcal{Y}$, $(3)$ implies that the image of $\varepsilon_I^{0}:I\to\ell(I)$ belongs to $\mathcal{Y}$. Further, since $\varepsilon_I^{0}:I\to \ell(I)$ is the unit adjunction of $I$, $\Img(\varepsilon_I^{0})=\ell(I)$. This shows $X^I=0$. Now, by Lemma \ref{EOff}, $\Db{i}$ is fully faithful since $\mathcal{A}$ has enough injectives. Thus $(1)$ holds.
 $\square$

\subsection{Nonsingular rings  and commutative noetherian rings\label{sect4.3}}

In this section we will construct derived decomposition by applying Proposition \ref{Localizing} to localizing subcategories of the modules over left nonsingular rings and commutative noetherian rings, respectively.

First, we consider left nonsingular rings (see \cite[Chapter 1]{Goodearl}).

Let $R$ be a ring and $M$ be an $R$-module with a submodule $N$. Recall that $M$ is an \emph{essential extension} of
$N$ (or $N$ is an \emph{essential submodule} of $M$) if every nonzero submodule of $M$ has nonzero intersection with $N$.
Recall that the injective envelope of $N$ is just an essential extension $M$ of $N$ with $M$ an injective module. As before, $M$ is denoted by $E(N)$. The set of all essential submodules of $_RR$ is denoted by $\mathscr{S}(R)$. A class $\mathcal{U}$ of $R$-modules is said to be \emph{closed under essential extensions} in $R\Modcat$ provided that $M\in\mathcal{U}$ whenever $M$ is an essential extension of a module $N\in\mathcal{U}$.

For an $R$-module $M$, we define
$Z(M):=\{x\in M\mid Ix=0\;\mbox{for some}\; I\in\mathscr{S}(R)\}$. This is a submodule of $M$ and called the \emph{singular submodule} of $M$.
The module $M$ is called a \emph{singular} module if $Z(M)=M$; and a \emph{nonsingular} module if $Z(M)=0$. The ring $R$ is said to be \emph{left nonsingular} if $_RR$ is a nonsingular module. Examples of left nonsingular rings include left semi-hereditary rings, direct products of integral domains, semiprime left Goldie rings and commutative semiprime rings (see \cite{Goodearl} for more examples).

\begin{Koro}\label{DADS}
Let $R$ be a ring, $\mathcal{X}$ be the full subcategory of singular modules in $R\Modcat$, and $\mathcal{Y}$ be the full subcategory of $R\Modcat$ consisting of all direct summands of arbitrary products of copies of $E(_RR)$. If $R$ is left nonsingular, then $\big(\mathcal{X}, \mathcal{Y}\big)$ is a derived decomposition of $R\Modcat$.
\end{Koro}

To show Corollary \ref{DADS}, we need the following basic properties of singular and nonsingular modules (see \cite[Propositions 1.20 and 1.22]{Goodearl} for proofs).

\begin{Lem} \label{singular}
Let $R$ be a ring, $\mathcal{X}$ be the full subcategory of singular modules in $R\Modcat$, and $M$ be an $R$-module.

$(1)$ $M\in \mathcal{X}$ if and only if $M$ is isomorphic to the quotient $X/Y$ of an essential extension $Y\subseteq X$.

$(2)$ $M\in \mathcal{X}^{\bot 0}$ if and only if $M$ is nonsingular.

$(3)$ $\mathcal{X}$ is closed under submodules, quotients and direct sums in $R\Modcat$; and $\mathcal{X}^{\bot 0}$
is closed under submodules, direct products, extensions and essential extensions in $R\Modcat$.
\end{Lem}

In general, the full subcategory $\mathcal{X}$ of singular $R$-modules may not be closed under extensions in $R\Modcat$. Nevertheless, the next lemma, taken from \cite[Propositions 1.23 and 2.12]{Goodearl}, provides a positive situation.

\begin{Lem}\label{nonsingular}
Let $R$ be a left nonsingular ring and $\mathcal{X}$ be the full subcategory of singular $R$-modules. Then

$(1)$ $\mathcal{X}$ is closed under extensions and essential extensions in $R\Modcat$.

$(2)$ $\mathcal{X}={^{\bot 0}}(\mathcal{X}^{\bot 0})={^{\bot 0}}E(R)$.

$(3)$ Let $M$ be an $R$-module. Then $M/Z(M)\in\mathcal{X}^{\bot 0}$. Moreover, $M$ is nonsingular if and only if
$M$ can be embedded in a direct product of copies of $E(R)$.
\end{Lem}

\begin{Rem}
If $R$ is a left nonsingular ring, then the subcategory $\mathcal{X}$ of singular $R$-modules is a localizing subcategory of $R\Modcat$ by Lemmas \ref{singular} and \ref{nonsingular}. Hence
$(\mathcal{X}, \mathcal{X}^{\bot 0})$ is a hereditary torsion pair in $R\Modcat$ (see \cite[Section 1, p.13]{BI} for definition). Moreover, $E(R)$ is nonsingular, while $E(R)/R$ is singular.
\end{Rem}

\noindent {\bf Proof of Corollary \ref{DADS}.}
Let  $R$ be a left nonsingular ring. We show $\mathcal{Y}=\mathcal{X}^{\bot 0}\cap \mathscr{I}(R\Modcat)=\mathcal{X}^{\bot}=\mathcal{X}^{\bot {0,1}}.$

Since $E(_RR)$ is injective and nonsingular, $\mathcal{Y}\subseteq \mathcal{X}^{\bot 0}\cap \mathscr{I}(R\Modcat)$.
The converse inclusion $\mathcal{X}^{\bot 0}\cap \mathscr{I}(R\Modcat)$ $\subseteq \mathcal{Y}$ follows from Lemmas \ref{singular}(2) and \ref{nonsingular}(3). Thus $\mathcal{Y}=\mathcal{X}^{\bot 0}\cap \mathscr{I}(R\Modcat)$.

Clearly, $\mathcal{Y}\subseteq\mathcal{X}^{\bot}\subseteq\mathcal{X}^{\bot {0,1}}$.
So, to show the other equalities, it is enough to show $\mathcal{X}^{\bot {0,1}}\subseteq \mathcal{Y}$.
To this purpose, we first prove a general result:

$(\ast\ast)\; $ If $0\to M\lraf{f} I\lraf{g} J$ is an exact sequence of $R$-module such that $I$ is injective and $J$ is nonsingular, then $M$ is injective.

In fact, since $f$ is injective and $_RI$ is injective, there is another injective $R$-module $K$ such that $\Coker(f)\simeq (E(M)/M)\oplus K$. Moreover, since $J$ is nonsingular, $\Img(g)$ is also nonsingular by Lemma \ref{singular}(3). It follows from $\Coker(f)=\Img(g)$ that $E(M)/M$ is nonsingular. However, $E(M)/M$ is
singular by Lemma \ref{singular}(1). This implies that $M=E(M)$ is injective.

\smallskip
By Lemma \ref{Compare}, the category $\mathcal{X}^{\bot {0,1}}$ consists of all $R$-modules $Y$ which has a minimal injective presentation
$0\to Y\to I_0\to I_1$ with $I_0, I_1\in\mathcal{Y}$. Since $I_0$ is injective and $I_1$ is nonsingular, it follows from $(\ast\ast)$ that $Y$ is injective. Thus $Y\in\add(I_0)\subseteq\mathcal{Y}$. This shows $\mathcal{X}^{\bot {0,1}}\subseteq \mathcal{Y}$.

As a consequence of $(\ast\ast)$, there is an isomorphism $I\simeq M\oplus \Img(g)$. Thus $\Img(g)$ is injective. This implies
$\Img(g)\in\add(J)$. Now, if $J\in\mathcal{Y}$, then $\Img(g)\in\mathcal{Y}$. Hence Corollary \ref{DADS} follows from Proposition \ref{Localizing}(3). $\square$

In the sequel, we deal with commutative noetherian rings.

Let $R$ be a commutative noetherian ring and $\Spec(R)$ be the prime spectrum of $R$.
For a multiplicative subset $\Sigma$ of $R$, we denote by $\Sigma^{-1}R$ the localization of $R$ at $\Sigma$. For $\mathfrak{p}\in\Spec(R)$, let $R_{\mathfrak p}$ be the localization of $R$ at the set
$R\setminus\mathfrak{p}$.
We always identify $\Spec(\Sigma^{-1}R)$ with the subset of all prime ideals $\mathfrak{p}$ of $R$ satisfying $\mathfrak{p}\cap \Sigma=\emptyset$. We also regard $(\Sigma^{-1}R)\Modcat$ as a full subcategory of $R\Modcat$ in the sense that an $R$-module $M$ belongs to $(\Sigma^{-1}R)\Modcat$ if and
only if $M\simeq (\Sigma^{-1}R)\otimes_R M$ as $R$-modules.
Let $\Ass(M)$ be the set of prime ideals $\mathfrak{p}$ of $R$ such that $R_{\mathfrak p}$ is isomorphic to a submodule of $M$. Thus $\Ass(M)=\Ass(E(M))$, where $E(M)$ is an injective envelope of $M$. The \emph{support} of $M$, denoted by $\Supp(M)$, is by definition the set of prime ideals $\mathfrak{p}$ of $R$ satisfying $\Tor_i^R(R_{\mathfrak p}/\mathfrak{p}R_{\mathfrak p},\, M)\neq 0$ for some $i\in\mathbb{N}$ (see \cite{Foxby}). In general, $\Ass(M)\subseteq \Supp(M)\subseteq \{\mathfrak{p}\in \Spec(R)\mid M_\mathfrak{p}\neq 0\}.$
The second inclusion is an equality if the module $M$ is finitely generated.
Note that $\Supp(M)$ is the union of the subsets $\Ass(I)$ of $\Spec(R)$, where $I$ runs over all those injective $R$-modules that appear
in a minimal injective resolution of $M$ (see \cite[Remark 2.9]{Foxby} or \cite[Lemma 3.3]{Krause}). In particular, if $M$ is injective, then $\Ass(M)=\Supp(M)$.

The following hold for a commutative noetherian ring $R$:

(a) Each injective $R$-module is a direct sum of indecomposable injective $R$-modules.

(b)
$\{E(R/\mathfrak{p})\mid \mathfrak{p}\in\Spec(R)\}$ is a complete set of non-isomorphic indecomposable injective $R$-modules.

(c) For $\mathfrak{p}, \mathfrak{q}\in\Spec(R)$, $\Hom_R(E(R/\mathfrak{p}), E(R/\mathfrak{q})\neq 0$ if and only if $\mathfrak{p}\subseteq\mathfrak{q}$
(see \cite[Theorems 3.3.7 and 3.3.8]{EJ}).

Let $\mathscr{S}$ be a full subcategory of $R\Modcat$ and let $\Phi$ be a subset of $\Spec(R)$. We define
$$\Supp(\mathscr{S}):=\bigcup_{M\in\mathscr{S}}\Supp(M) \quad\mbox{and}\quad  \Supp^{-1}(\Phi):=\{M\in R\Modcat\mid \Supp(M)\subseteq \Phi\}.$$
Gabriel's classification of localizing subcategories (see \cite[p. 425]{Gabriel}) conveys that the map  $\Supp$ induces a bijection between the set of localizing subcategories of $R\Modcat$ and the set of specialization closed subsets of $\Spec(R)$. The inverse of $\Supp$ is just given by $\Supp^{-1}$.
This was extended in \cite[Theorem 3.1]{Krause} to a bijection (with the same maps) between the set of abelian full subcategories of $R\Modcat$ closed under extensions and arbitrary direct sums, and the set of coherent subsets of $\Spec(R)$.

A subset $\Phi$ of $\Spec(R)$ is said to be \emph{specialization closed} provided that if $\mathfrak{p}, \mathfrak{q}\in\Spec(R)$ and
$\mathfrak{p}\subseteq\mathfrak{q}$, then $\mathfrak{p}\in\Phi$ implies $\mathfrak{q}\in\Phi$; \emph{coherent} provided that each homomorphism $I^0\to I^1$ between injective $R$-modules with $\Ass(I^0)\cup \Ass(I^1)\subseteq\Phi$ can be completed to an exact sequence $I^0\to I^1\to I^2$ such that $I^2$ is injective and $\Ass(I^2)\subseteq\Phi$ (see \cite[Section 3]{Krause}). Examples of coherent subsets are specialization closed subsets and $\Spec(\Sigma^{-1}R)$.
For futher information on coherent subsets, we refer the reader to \cite[Section 4]{Krause}.

An application of Proposition \ref{Localizing} is the following

\begin{Koro}\label{Localizing commutative}
Let $R$ be a commutative noetherian ring, $\Phi$ be a specialization closed subset of $\Spec(R)$, and  $\Phi^{c}:=\Spec(R)\setminus\Phi$.
Then the pair $\big(\Supp^{-1}(\Phi), \Supp^{-1}(\Phi^c)\big)$ is a derived decomposition of $R\Modcat$ if and only if $\Phi^c$ is coherent.
\end{Koro}

{\it Proof.}
Let $\mathcal{X}:=\Supp^{-1}(\Phi)$ and $\Ass^{-1}(\Phi^c):=\{M\in R\Modcat\mid \Ass(M)\subseteq \Phi^c\}$. We first show that
$\mathcal{X}^{\bot 0}=\Ass^{-1}(\Phi^c)$.

Let $U\in\mathcal{X}$ and $V\in\Ass^{-1}(\Phi^c)$. Then $\Supp(E(U))\subseteq \Supp(U)\subseteq\Phi$ and $\Ass(E(V))=\Ass(V)\subseteq \Phi^c$.
If $\Hom_R(E(U), E(V))\neq 0$, then there is a non-zero homomorphism from a direct summand $E(R/\mathfrak{p})$
of $E(U)$ to a direct summand $E(R/\mathfrak{q})$ of $E(V)$, where $\mathfrak{p}\in\Phi$ and $\mathfrak{q}\in\Phi^c$.
In this case, we have $\mathfrak{p}\subseteq\mathfrak{q}$. This is contradictory to the assumption that $\Phi$ is specialization closed.
Thus $\Hom_R(E(U), E(V))= 0$. This implies $\Hom_R(U,V)=0$ and shows $\Ass^{-1}(\Phi^c)\subseteq \mathcal{X}^{\bot 0}$.
To verify $\mathcal{X}^{\bot 0}\subseteq\Ass^{-1}(\Phi^c)$, we take $W\in \mathcal{X}^{\bot 0}$ and $\mathfrak{a}\in \Ass(W)$.
Then  $R/\mathfrak{a}$ is isomorphic to a nonzero submodule of $W$. If $\mathfrak{a}\in\Phi$, then $\Supp(R/\mathfrak{a})=\{\mathfrak{b}\in\Spec(R)\mid \mathfrak{a}\subseteq\mathfrak{b}\}\subseteq \Phi$ since $\Phi$ is specialization closed.
This implies $R/\mathfrak{a}\in \mathcal{X}$, and therefore $\Hom_R(R/\mathfrak{a}, W)=0$. This is a contradiction. Thus $\mathfrak{a}\in\Phi^c$.

Let $\mathcal{Y}:=\mathcal{X}^{\bot}$ and $\mathcal{E}:=\mathcal{X}^{\bot 0}\cap \mathscr{I}(R\Modcat)$. Then $\mathcal{E}=\Ass^{-1}(\Phi^c)\cap \mathscr{I}(R\Modcat)=\mathcal{Y}\cap \mathscr{I}(R\Modcat)$. By Lemma \ref{Compare}, $\mathcal{Y}$ consists of all $R$-modules $M$ which has a minimal injective resolution $0\to M\to I_0\to I_1\to I_2\to\cdots $ such that $I_i\in\mathcal{E}$ for all $i\geq 0$. Since $\Supp(M)=\bigcup_{i\geq 0}\Ass(I_i)$, $\mathcal{Y}=\Supp^{-1}(\Phi^c)$. Observe that, for any $\Psi\subseteq\Spec(R)$, $\Supp^{-1}(\Psi)$ is always closed under direct sums in $R\Modcat$ because $\Supp(\bigoplus_{j\in J}M_j)=\bigcup_{j\in J}\Supp(M_j)$ for any family $\{M_j\}_{j\in J}$ of $R$-modules with $J$ an index set. In particular, $\mathcal{Y}$ is closed under direct sums in $R\Modcat$.

If $(\mathcal{X}, \mathcal{Y})$ is a derived decomposition of $R\Modcat$, then
$\mathcal{Y}$ is an abelian full subcategory of $R\Modcat$ and closed under both extensions and direct sums.
In this case, $\Phi^c$ is coherent. This shows the necessity in Corollary \ref{Localizing commutative}.

Conversely, suppose that $\Phi^c$ is coherent. Let $f:I^0\to I^1$ be a homomorphism between injective $R$-modules with $I^1\in\mathcal{Y}$.
By Proposition \ref{Localizing}(2), we need to extend $f$ to an exact sequence $I^0\to I^1\to I^2$ in $R\Modcat$ with $I^2\in\mathcal{E}$.
This can be done if $I^0\in\mathcal{E}$ since $\Phi^c$ is coherent. For the general case, we decompose $I^0$ into a direct sum of
indecomposable injective modules. Recall that $\{E(R/\mathfrak{p})\mid \mathfrak{p}\in\Spec(R)\}$ is a complete set of
isomorphism classes of indecomposable injective $R$-modules and that $\Ass(E(R/\mathfrak{p}))=\Supp(E(R/\mathfrak{p}))=\{\mathfrak{p}\}$. Consequently,
$E(R/\mathfrak{p})$ belongs to either $\mathcal{X}$ or $\mathcal{Y}$. This yields a decomposition $I^0=X\oplus Y$ with $X\in\mathcal{X}$ and
$Y\in\mathcal{Y}$. Since $\Hom_R(X,I^1)=0$, $f=(0, g)$, where $g:Y\to I^1$ is the restriction of $f$ to $Y$.
Clearly, $g$ is a homomorphism between modules in $\mathcal{E}$. Now, we first extend $g$ and then $f$
to an exact sequence $I^0\to I^1\to I^2$ in $R\Modcat$ with $I^2\in\mathcal{E}$. Thus $(\mathcal{X}, \mathcal{Y})$ is a derived decomposition of $R\Modcat$. $\square$

\medskip
In Corollary \ref{Localizing commutative}, when $\Phi^c$ is coherent, $\mathcal{Y}:=\Supp^{-1}(\Phi^c)$ is an abelian full subcategory of $R\Modcat$
closed under direct sums, and the inclusion $j:\mathcal{Y}\to R\Modcat$ has a left adjoint $\ell: R\Modcat\to\mathcal{Y}$.
Let $S=\End_R(\ell(R))$ and let $\lambda:R\to S$ be the ring homomorphism induced from the functor $\ell$.
Thanks to \cite[Proposition 3.8]{GL}, $\lambda$ is a ring epimorphism, inducing an equivalence of abelian categories: $S\Modcat \lraf{\simeq} \mathcal{Y}$, and
$S$ is a flat $R$-module since $\ell$ is exact. Thus $\lambda$ is a flat ring epimorphism (see also \cite{AMTTV} for further details). Consequently,
$S$ is also a commutative noetherian ring. So we can apply Corollary \ref{Localizing commutative} (for example, via localizations) to $S$ and obtain
a derived decomposition of $S\Modcat$. By iterating this procedure, we can stratify $R\Modcat$ as a sequence of
derived decompositions of the module categories over commutative rings. In particular, when the Krull dimension of $R$ is at most $1$, a derived stratification (see Definition \ref{stratification}) of $R\Modcat$ can be constructed explicitly.

The following is a restatement of Corollary \ref{cor1.4}.

\begin{Koro}\label{special cases}
Let $R$ be a commutative noetherian ring.

$(1)$ Suppose that $\Phi$ is a specialization closed subset of $\Spec(R)$. If the Krull dimension of $R$ is at most $1$, then $\big(\Supp^{-1}(\Phi), \Supp^{-1}(\Phi^c)\big)$ is a derived decomposition of $R\Modcat$, where $\Phi^{c}:=\Spec(R)\setminus\Phi$.

$(2)$ Let $\Sigma$  be a multiplicative subset of $R$  and $\Phi:=\{\mathfrak{p}\in\Spec(R)\mid\mathfrak{p}\cap\Sigma\neq \emptyset\}.$ Then
$\big(\Supp^{-1}(\Phi), (\Sigma^-R)\Modcat\big)$ is a derived decomposition of $R\Modcat$.
\end{Koro}

{\it Proof.} By \cite[Theorem 1.2]{Krause}, if the Krull dimension of $R$ is at most $1$, then every subset of $\Spec(R)$ is coherent.
Thus $(1)$ follows from Corollary \ref{Localizing commutative}. For $(2)$, we recall that $\Spec(\Sigma^{-1}R)$ is identified with the subset of all prime ideals $\mathfrak{p}$ of $R$ satisfying $\mathfrak{p}\cap \Sigma=\emptyset$. Thus $\Phi^c=\Spec(\Sigma^{-1}R)$ and is coherent.
Further, $\Supp^{-1}(\Phi^c)=(\Sigma^-R)\Modcat$. Now $(2)$ is a consequence of Corollary \ref{Localizing commutative}.
Remark that $(2)$ also follows from Corollary \ref{commtative}. $\square$

\begin{Koro}\label{decomposition}
Suppose that $R$ is a commutative noetherian ring of Krull dimension at most $1$. Let ${\rm Max}(R)$ be the set of maximal ideals of $R$ and
let ${\rm Min}(R)$ be the set of prime ideals of $R$ which are not maximal. Then

$(1)$ $\big(\Supp^{-1}({\rm Max}(R)),\, \Supp^{-1}({\rm Min}(R))\big)$ is a derived decomposition of $R\Modcat$.

$(2)$ There are equivalences of abelian categories:
$$\Supp^{-1}({\rm Max}(R))\lraf{\simeq}\prod_{\mathfrak{m}\in{\rm Max}(R)}\Supp^{-1}(\{\mathfrak{m}\})\quad \mbox{and}\quad
\Supp^{-1}({\rm Min}(R))\lraf{\simeq}\prod_{\mathfrak{p}\in{\rm Min}(R)}R_\mathfrak{p}\Modcat,$$
where $\prod$ denotes the direct product of abelian categories.

$(3)$ Both $\Supp^{-1}(\{\mathfrak{m}\})$ and $R_\mathfrak{p}\Modcat$ are abelian simple for any $\mathfrak{m}\in{\rm Max}(R)$ and $\mathfrak{p}\in{\rm Min}(R)$.
\end{Koro}

{\it Proof.} $(1)$ Clearly, ${\rm Max}(R)$ is a specialization closed subset of $\Spec(R)$. Since the Krull dimension of $R$ is at most $1$, the statement $(1)$
follows from Corollary \ref{special cases}(1).

$(2)$ To show the equivalences, we first establish the following

\smallskip
{\bf General fact.} Let $\Phi$ be a subset of $\Spec(R)$ with the property $(\lozenge)$: for $\mathfrak{p},\mathfrak{q}\in\Phi$,  $\mathfrak{p}\subseteq\mathfrak{q}$
implies  $\mathfrak{p}=\mathfrak{q}$.  Then there is an equivalence of abelian categories:
$\prod_{\mathfrak{p}\in\Phi}\Supp^{-1}(\{\mathfrak{p}\})\lraf{\simeq}\Supp^{-1}(\Phi).$

\smallskip
In fact, given an $R$-module $M$ with  a minimal injective resolution $0\to M\to I_0\to I_1\to I_2\to \cdots $,
we always have $\Supp(M)=\bigcup_{i\geq 0}\Ass(I_i)=\bigcup_{i\geq 0}\Supp(I_i)$. Consequently, $M\in\Supp^{-1}(\Phi)$ if and only if
$I_i\in \Supp^{-1}(\Phi)$ for all $i\geq 0$. Recall that, for each $\mathfrak{a}\in\Spec(R)$, if $\mathfrak{a}\nsubseteq\mathfrak{b}\in\Spec(R)$, then $\Hom_R(E(R/\mathfrak{a}), E(R/\mathfrak{b})=0$. Moreover, $R_\mathfrak{a}\Modcat$ is regarded as an abelian full subcategory of $R\Modcat$
and $E(R/\mathfrak{a})\in R_\mathfrak{a}\Modcat$. So, the property $(\lozenge)$ implies

$(1)$ $M\in\Supp^{-1}(\Phi)$ if and only if $M\simeq\bigoplus_{\mathfrak{p}\in\Phi}M_\mathfrak{p}$ with $M_\mathfrak{p}\in\Supp^{-1}(\{\mathfrak{p}\})$,
where $M_\mathfrak{p}$ stands for the localization of $M$ at $\mathfrak{p}$; and

$(2)$ If $\mathfrak{p},\mathfrak{q}\in\Phi$ and $\mathfrak{p}\neq \mathfrak{q}$, then $\Hom_R(X,Y)=0$ for all $X\in \Supp^{-1}(\{\mathfrak{p}\})$ and $Y\in \Supp^{-1}(\{\mathfrak{q}\})$.

By $(1)$ and $(2)$, one can check that
the functor $\bigoplus: \prod_{\mathfrak{p}\in\Phi}\Supp^{-1}(\{\mathfrak{p}\})\to\Supp^{-1}(\Phi)$, given by taking direct sums in $R\Modcat$,
is an equivalence of abelian categories.

\smallskip
Since the Krull dimension of $R$ is at most $1$, both ${\rm Max}(R)$ and ${\rm Min}(R)$ have the property $(\lozenge)$.
Moreover, if $\mathfrak{p}$ is a minimal prime ideal of $R$, then it follows from  $\Supp(M)=\bigcup_{i\geq 0}\Ass(I_i)$ that $\Supp^{-1}(\{\mathfrak{p}\})=R_\mathfrak{p}\Modcat$. Note that ${\rm Min}(R)$ consists of all minimal prime ideals of $R$
which are not maximal. Now, the existence of equivalences in Corollary \ref{decomposition} follows from the general fact.

$(3)$ If $(\mathcal{X},\mathcal{Y})$ is an Ext-orthogonal decomposition of an abelian category $\mathcal{A}$ with arbitrary direct sums, then $\mathcal{X}= {}^{\bot}\mathcal{Y}$ and $\mathcal{X}$ is closed under both extensions and arbitrary direct sums (see the dual of \cite[Theorem 1, p. 116]{Maclane}). It follows from \cite[Theorem 3.1]{Krause} that, for any $\mathfrak{a}\in\Spec(R)$, the abelian category $\Supp^{-1}(\{\mathfrak{a}\})$ does not contain non-trivial
abelian full subcategory which is closed under extensions and arbitrary direct sums.  This implies that $\Supp^{-1}(\{\mathfrak{a}\})$ is abelian simple. Thus $(3)$ holds.
$\square$

\begin{Rem}
The Krull dimension in Corollary \ref{decomposition} cannot be relaxed. For example, if $R$ is the algebra of formal power series over a field $k$ in two variables $x$ and $y$. Then $R$ is a local ring of Krull dimension $2$ and the maximal ideal $\mathfrak{m}$ of $R$ generated by $x$ and $y$ is specialization closed, but
$\Spec(R)\setminus\{\mathfrak{m}\}$ is not coherent (see \cite[Example 3.2]{Krause}). By Corollary \ref{Localizing commutative},
$\big(\Supp^{-1}(\{\mathfrak{m}\}), \Supp^{-1}(\Spec(R)\setminus\{\mathfrak{m}\})\big)$ is not a derived decomposition of $R\Modcat$.
\end{Rem}

{\footnotesize
Hongxing Chen,

School of Mathematical Sciences, Capital Normal University, 100048
Beijing, CHINA

{\tt Email: chx19830818@163.com}

\bigskip
Changchang Xi,

School of Mathematical Sciences, Capital Normal University, 100048
Beijing, CHINA; and

School of Mathematics and Information Science, Henan Normal University, Henan, CHINA

{\tt Email: xicc@cnu.edu.cn}}

\end{document}